\documentclass[11pt,leqno]{article}
\usepackage{amsmath,txfonts,fourier}
\usepackage{graphics,hyperref}
\numberwithin{equation}{section}
\usepackage{color}        
\newtheorem{theorem}{Theorem}[section]
\newtheorem{remark}[theorem]{Remark}
\newtheorem{lemma}[theorem]{Lemma}
\newtheorem{proposition}[theorem]{Proposition}
\newtheorem{corollary}[theorem]{Corollary}
\newtheorem{definition}[theorem]{Definition}
\newcommand{\bos}{\begin{remark}\rm}
\newcommand{\eos}{\end{remark}}
\newcommand{\bte}{\begin{theorem}}
\newcommand{\ete}{\end{theorem}}
\newcommand{\bpr}{\begin{proposition}}
\newcommand{\epr}{\end{proposition}}
\newcommand{\bdf}{\begin{definition}\rm}
\newcommand{\edf}{\end{definition}}
\newcommand{\bco}{\begin{corollary}}
\newcommand{\eco}{\end{corollary}}
\newcommand{\ble}{\begin{lemma}}
\newcommand{\ele}{\end{lemma}}
\newcommand{\bdm}{\begin{displaymath}}
\newcommand{\edm}{\end{displaymath}}
\newcommand{\beq}{\begin{equation}}
\newcommand{\eeq}{\end{equation}}
\newcommand{\bary}{\begin{eqnarray}}
\newcommand{\eary}{\end{eqnarray}}
\newcommand{\barr}{\begin{eqnarray*}}
\newcommand{\earr}{\end{eqnarray*}}
\newcommand{\bdim}{{\noindent{\bf Proof.}\quad}}
\newcommand{\edim}{{\unskip\nobreak\hfil\penalty50
\hskip2em\hbox{}\nobreak\hfil\mbox{\rule{1ex}{1ex} \qquad}
\parfillskip=0pt \finalhyphendemerits=0\par\medskip}}

\newcommand{\dys}{\displaystyle}
\newcommand{\dyle}{\displaystyle}

\newcommand{\rife}[1]{(\ref{#1})}

\def\R{\varmathbb{R}}
\def\RN{\R^{N}}

\def\hsob{H^{1}(\RN)}

\newcommand{\intr}{\int_{\RN}}
\newcommand{\grad}{\nabla}
\newcommand{\pa}{\partial}

\newcommand{\diff}{\!\setminus\!}

\newcommand{\jbe}{J_{\eps}}

\newcommand{\ov}[1]{\overline{#1}}

\newcommand{\vep}{\varepsilon}
\newcommand{\eps}{\varepsilon}

\newcommand{\vfi}{\varphi}

\renewcommand{\H}{\varmathbb{H}^{1}}

\newcommand{\opic}{{\rm o}}

\begin{document}

\title{Singularly perturbed elliptic  problems with
 nonautonomous asymptotically linear nonlinearities}

\author{Liliane A. Maia\thanks{Research partially 
supported by CNPq/Brazil and
FAPDF/Brazil},\; Eugenio Montefusco$^{\dag}$,\; 
Benedetta Pellacci\thanks{Research partially 
supported by FIRB project 2012 {\it Dispersive Dynamics} and 
PRIN project 2009-WRJ3W7}}
\date{}
\maketitle

\begin{abstract}
\noindent
We consider a class of singularly perturbed elliptic  
problems with nonautonomous asymptotically linear nonlinearities.
The dependence on the spatial coordinates co\-mes from the presence
of a potential and of a function representing a saturation effect.
We investigate the existence of nontrivial nonnegative solutions
concentrating around local minima of both the potential and of the
saturation function. 
Necessary conditions to locate the possible concentration points
are also given.
\end{abstract}

\vspace{.5cm}

\section{Introduction}
In this paper we study the existence of positive 
solutions of the problem
\beq\label{Eq:eps}
\begin{cases}
\dys -\eps^{2}\Delta u+V(x)u^{2}
=\frac{u^{3}}{1+s(x)u^{2}} &\text{in } \R^{N},
\\
\dys \,u\,\in H^{1}(\R^{N}),
\tag{$P_{\eps}$}\end{cases}
\eeq
for $N\geq2$, $\eps>0$ a small parameter and $V,\,s:\R^{N}\to \R$  
H\"older continuous functions such that
\beq\label{s0}
s(x)\geq\alpha>0 \qquad \forall x\in\RN,
\eeq
\beq\label{V0}
V(x)\geq \mu>0, \qquad \forall x\in\RN.
\eeq
It is well known that every positive solutions $u_{\eps}$ of 
\eqref{Eq:eps} generates a standing wave, i.e. 
$\phi_{\eps}(x,t)=u_{\eps}(x)e^{-iEt/\hbar}$ solution of
\beq\label{schintr}
i\hbar \partial_{t}\phi+\frac{\hbar^{2}}{2m}\Delta \phi-W(x)\phi
=\frac{\phi^{3}}{1+s(x)\phi^{2}}
\eeq
for $W=V+
E$, $\eps^{2}=\hbar/2m$.
\\
Problem \eqref{schintr} represents the propagation of a light pulse
along a saturable medium.  A typical class of saturable medium
is constituted by the photorefractive crystals,
one of the most preferable materials to observe the propagation of a 
light beam, because of their slow response to the propagation, making 
easier the observation. When a beam passes through  these materials
its refractive index changes so that the light remains confined 
and solitons are generated.
When observing light propagation through these media one 
can see a {\it saturation effect}: it is possible to increase
the amplitude of the generated solitons by increasing
light intensity up to a critical bound characteristic of the material.
This kind of interaction is not well represented by 
the usual Schr\"odinger
equation, so that this model is replaced by \eqref{schintr}
where the usual autointeraction represented by the cubic power
is prevalent for "small" $u$, while a linear interaction,
$u/s(x)$, is predominant for "large" $u$.
Moreover, aiming to analyze the observation through
different materials we admit a possible change 
of the saturation feature in dependence on the spatial
coordinates, which may happen observing the propagation
along different material.

An interesting and largely studied class of 
solutions of \eqref{Eq:eps} is the family of 
semiclassical states, that are families $u_{\eps}$ 
with a spike shape concentrating around some points of $\R^{N}$
for $\eps$ sufficiently small.
There is a broad variety of contributions concerning the existence 
of this kind of solutions for the equation
\beq\label{ell}
-\eps^{2}\Delta u+V(x)u=f(x,u).
\eeq
For $f(x,t)=t^{3}$, 
the first contribution on the subject in the one dimensional case
is due to Floer and Weinstein \cite{flwe} who show the existence
of a solution $u_{\eps}$  concentrating
around any given $x_{0}$ nondegenerate  critical point of 
$V(x)$. Their result has been extended in higher dimension
in \cite{oh1,oh2} for $f(x,t)=|t|^{p-1}t$ with $1<p<(N+2)/(N-2)$.
The common approach used in this papers is a 
Lyapunov-Schmidt reduction, consisting in a local bifurcation
tuype result, which relies on the uniqueness and 
nondegeneracy of the ground state solution of the autonomous 
problem
\beq\label{autointr}
-\Delta v+V(x_{0})v=f(x_0,v).
\eeq
The Lyapunov-Schmidt procedure or more general finite-dimensional
reductions methods have been used to find solutions concentrating
around any \(x_0\) isolated minimum (or maximum) point
with possibly polynomial degeneration of $V$ in
\cite{ambaci}, and then around different stable critical points 
(see \cite{li,gr,pi,amfema} and the references therein).

A different approach to this is to find a solution 
$u_{\eps}$ for $\eps$ positive and then study its 
asymptotic behavior for $\eps$ tending to zero. 
This procedure has been firstly used by Rabinowitz in \cite{ra}
assuming that  $\inf V(x)<\liminf_{|x|\to +\infty}V(x)$ 
and proving concentration around a local minimum point of $V$.
This philosophy has been improved in \cite{defe,defe2}, where
it is shown, by means of a penaliztion argument, that it is
sufficient to assure a local condition on the potential:
there exists a bounded open set $\Lambda$ such that
$$
\inf_{\Lambda}V<\inf_{\partial \Lambda}V.
$$
As for the reduction method also this procedure has been used 
to extend the existence and concentration result in many 
different directions (see \cite{defe2, beje, bova}).
\\
When passing in \eqref{autointr} from $f(t)=k(x) t^{3}$ to 
$f(x,t)=t^{3}/(1+s(x)t^{2})$ many differences arises.
First of all, thanks to \eqref{s0}, we do not have a critical 
exponent as  $|f(t)|<t/\alpha$. Moreover, as $f$ is 
asymptotically linear, the action functional $I_\eps$, 
defined in 
\beq\label{defh}
\H=\left\{u\in H^{1}(\R^{N})\,:\, V(x)u^{2}\in L^{1}(\R^{N})\right\},
\eeq
by
\[
I_{\eps}(u):=\frac{1}2\intr \left[\eps^{2} |\nabla u(x)|^{2}dx+
V(x)u^{2}(x)\right]dx-\intr F(x,u(x))dx,
\]
for $F(x,t)$ given by
\beq\label{defF}
F(x,t)= \frac1{2s(x)} t^2 -\frac{1}{2s^{2}(x)}\ln(1+s(x)t^2),
\eeq
may present different geometric behavior in dependence 
of $V$ and $s$, e.g. if $V(x)s(x)>1$ for every $x\in \R^{N}$, 
$I_\eps$ is always positive, convex and has only a global 
minimum at $u\equiv 0$.  
For $V$ and $s$ constant and such that $Vs<1$ in \cite{stzo} 
it is proved the existence of a positive radially symmetric
solution which is showed to be unique according to
\cite{seta,sezo}. 
Regarding the existence of semiclassical states, in \cite{jeta2} 
it is studied this kind of problem for general autonomous 
nonlinearity $f(x,t)=f(t)$, asymptotically linear or not,
and it is shown the existence of a positive solution 
$u_{\eps}$ concentrating around a local minimum of $V$ via 
variational methods and penalization arguments.
Here, being interested in the possible interaction between 
$V(x)$ and $s(x)$, we will deal with the following autonomous, 
or frozen, problem
\begin{equation}\label{limity}
\tag{$S_y$}
\begin{cases}
\dys-\Delta u + V(y)u = \frac{u^3}{1+s(y)u^{2}} & \text{in $\R^N$},
\\
u\in H^{1}(\R^{N}),
\end{cases}
\end{equation}
which has a solution if and only if $y$ belongs to the open 
set (see \cite{stzo})
\begin{equation}\label{defo}
\Omega=\{y\in \R^{N}\,:\, V(y)s(y)<1\}.
\end{equation}

Therefore, the set of possible concentration points  
is restricted from the beginning. 
As a further consequence, it is not possible to project
every $u\in H^{1}(\RN)$ on the Nehari manifold
\beq\label{ny}
{\mathcal N}_{y}:=\left\{u\in H^{1}\setminus\{0\},\,:\, 
\langle I'_{y}(u),u\rangle=0\right\},
\eeq
where $I_{y}$ is the autonomous, or frozen, functional
\beq\label{Iy}
I_{y}(v):=\frac12\intr |\nabla u|^{2}dx+
\frac1{2}V(y)\intr u^{2}dx
-\intr F(y,u(x))dx.
\eeq
Nevertheless, the nonlinearity $f(t)=t^{3}/(1+st^{2})$ is 
such that $f(t)/t$ is increasing w.r.t. $t$, ensuring the
uniqueness of the projecton whenever it exists.
\\
Another effect of the asymptotically linearity property
of $f$ is the loss of the well known Ambrosetti-Rabinowitz 
condition
$$
\exists\, \theta>2,\;\text{such that }\, \theta F(x,t)\leq f(x,t)t.
$$
This condition is useful in proving the boundedness of 
a Palais-Smale sequence. Here, we overcome this difficulty 
noticing that $f$ satisfies the so-called nonquadraticity 
condition
\beq\label{nonq}
f(x,t)t-2F(x,t)\geq 0,\;\text{and}\;
\lim_{|t|\to +\infty}f(x,t)t-2F(x,t)=+\infty.
\eeq
where $F$ is the primitive of $f$ (w.r.t. $x$) 
such that $F(x,0)=0$ (see Lemma \ref{NQpernoi}).

This condition enables us to show the boundedness of a Cerami
sequence. Our main result concerning sufficient conditions 
is stated in the following result, where we denote with
$B(z,r)$ the open ball centered at $z$ with radius $r$.

\bte\label{mainvs}
Assume condition \eqref{s0}, \eqref{V0}. Moreover, suppose that
there exists $z\in \Omega$ and $r>0$ such that either
\beq\label{minivs}
\left\{\begin{array}{c}
\dyle V_{0}=V(z)=\min_{B(z,r)}V (x)\leq\min_{\partial B(z,r)}V (x), 
\quad
s_{0}=s(z)=\min_{B(z,r)}s (x)<\min_{\partial B(z,r)}s (x),\\ 
\hbox{or} \\
\dyle V_{0}=V(z)=\min_{B(z,r)}V (x)<\min_{\partial B(z,r)}V (x), 
\quad
s_{0}=s(z)=\min_{B(z,r)}s (x)\leq\min_{\partial B(z,r)}s (x).
\end{array}\right.
\eeq
Then there exists $\eps_{0}>0$ such that, for every 
$0<\eps<\eps_{0}$, problem $(P_{\eps})$ admits a nontrivial  
solution $u_{\eps}\in\H$, $u_{\eps}\geq0$, such that 
the following facts hold:
\vskip3pt
\noindent(i) $u_{\eps}$ admits exactly one global maximum
point $x_{\eps}\in B(z,r)$;
\vskip3pt
\noindent(ii) 
$\dys \lim_{\eps\to 0} V(x_{\eps})=V_0$\; and\; 
$\dys \lim_{\eps\to 0} s(x_{\eps})=s_0;$ 
\vskip3pt
\noindent(iii) 
there exist $\mu_1,\mu_2>0$ such that, for every 
$x\in\R^N$,
$$
u_\eps(x)\leq \mu_1e^{-\mu_2\textstyle\frac{|x-x_\eps|}{\eps}}.
$$
\ete

Notice that, differently from the most studied case
$f(x,t)=k(x)|t|^{p-1}t$ (see \cite{wangzeng}), 
here concentration is produced around minimum points of 
$V$ and $S$, moreover notice that the strict inequality
is needed only on $s$ or $V$ not on both. 
Then, for example, one between $V$ or $s$ can be constant.
We can also prove an abstract concentration 
result around minimum points of the function
$\Sigma:\R^{N}\to \R$  defined by
\beq\label{defsigmagen}
\Sigma(y)=\begin{cases}
\dys\inf_{{\mathcal N}_{y}}I_{y} & y\in \Omega, \\
+\infty & y\not\in \Omega.
\end{cases}
\eeq

Thanks to the uniqueness of the ground state solution
of the autonomous problem \eqref{autointr} 
$\Sigma$ is regular in $\Omega$, nevertheless we cannot
obtain an explicit formula for $\Sigma$ because of the lack
of homogeneity of the autonomous problem.
More precisely, when $f(x,t)=K(x)|t|^{p-1}t$ one can derive
every solution of the equation \eqref{ell}
via a change of scale, from the unique positive 
solution of the equation $-\Delta u+u=|u|^{p-1}u$.
Here, this procedure cannot work as there are no nontrivial
solutions of the problem $-\Delta u+u=u^{3}/(1+u^{2})$,
so that we have no hope to find an explicit function
of $V$ and $s$, the critical points of which consituting 
the concentration set.
In Theorem \ref{mainnece} it is given a  necessary 
condition for the concentration to occur. 
Studying $\Sigma$ we realize that if $z$ is a 
concentration point, then
the gradient of $V$ and $s$ must be linearly dependent as it results 
in \cite{sesq} for a different class of problems.
Moreover, in our case the gradients $\nabla V(z)$ and 
$\nabla s(z)$ point in opposite directions and either 
$z$ is a common zero of $\nabla V(z)$ and $\nabla s(z)$
or $\nabla V(z),\nabla s(z)$ are both different from zero
and still $\nabla \Sigma (z)=0$.

\section{Setting of the Problem and Main Results}\label{setting}

In order to study $(P_{\eps})$ it is natural to 
introduce the Hilbert space $\H$, defined in \eqref{defh},
with norm $\|u\|_{\H}^2=\|u\|_{\eps,V}^{2}$, given by
$$
\|u\|_{\eps,V}^{2}=\eps^{2}
\|\nabla u\|_{2}^{2}+\int_{\R^{N}}V(x)u^{2}dx,
$$
where we denote with $\|\cdot\|_p$
the standard norm in $L^p=L^p(\R^{N})$ for $1\leq p\leq\infty$.
Thanks to condition \eqref{s0} we can say that the 
solutions of problem
\eqref{Eq:eps} correspond to the critical points of the $C^{1}$
functional  $I_{\eps}:H^{1}\to\R$ defined by
\beq\label{defI}
I_{\eps}(u)=\frac12\|u\|_{\eps,V}^{2}-\int_{\R^{N}}F(x,u(x))dx,
\eeq
for $F(x,t)$ given in \eqref{defF}.

It is easily  checked that $I_{\eps}$ is well defined 
and of class $C^{1}$ on $\H$.
A nontrivial solution of problem $(P_\vep)$ is a  $u_\vep\neq 0$
in $\H$, critical point of $I_{\eps}$.

For every $y\in \Omega$ (see \eqref{defo}) we can deduce from 
\cite{stzo, seta, sezo}  that there exists a unique, 
positive, radially symmetric least energy solution, 
denoted by $Q_{y}$, of the autonomous problem
frozen in $y$ \eqref{limity}. $Q_{y}$ is a critical point of 
the autonomous functional $I_{y}$, defined in \eqref{Iy}
and, denoting with $f(y,t)=\partial_{t}F(y,t)$, notice 
that $f(y,t)/t$ is an increasing function with respect to $t$. 
This monotonicity property is crucial in proving that the 
Mountain Pass level equals the minimum on the Nehari 
manifold (see Proposition 3.11 in \cite{ra}).
This equivalence will be often used in the sequel.

The first sufficient condition for the concentration 
Notice that, since for every continuous function $k(x)$
$$
\inf_{B(z,r)} k(x) \leq \min_{\pa B(z,r)} k(x) 
$$
in the inequality concerning $V$ in the first alternative
in \eqref{minivs}, it is assumed that the infimum of $V$
in $B(z,r)$ is actually achived in $z$ to occur is contained 
in Theorem \ref{mainvs}. With this respect the following 
comments are in order.

\bos
In \eqref{minivs} it is supposed that the minimum 
values $s_{0}$ and $V_{0}$ are achieved at $z$, 
the center of the ball. 
This can always be assumed without loss of generality, 
indeed condition \eqref{minivs} implies 
that $s_{0}$ and $V_{0}$ are achieved at a point 
$z_{1}\in B(z,r)\cap \Omega$, as $s(z_{1})V(z_{1})\leq s(z)V(z)<1$. 
Therefore, if $z_{1}\neq z$ we can replace $z$ with $z_{1}$ 
in the statement of the Theorem and in all the changes 
of variable in Section \ref{proof}, obtaining 
concentration around $z_{1}$.
On the other hand, it would be interesting to study the case 
in which the concentration occurs in different critical points 
of $V$ and $s$.
\eos

\bos
Notice that, since for every continuous function $k(x)$
$$
\inf_{B(z,r)} k(x) \leq \min_{\pa B(z,r)} k(x) 
$$
in the inequality concerning $V$ in the first alternative
in \eqref{minivs}, it is assumed that the infimum of $V$
in $B(z,r)$ is actually achived in $z$ and it may be
equal to the minimum of $V$ on the boundary (analogous
considerations hold for $s$ in the second alternativ in \eqref{mainvs}).
\eos

We can also prove the following general abstract result.
\bte\label{mainsigma}
Assume condition \eqref{s0}, \eqref{V0}. Moreover, suppose that
there exists $z\in \Omega$ and $r>0$ such that
\beq\label{minisigma}
\Sigma_{0}=\Sigma(z)=\min_{B(z,r)}\Sigma (x)
<\min_{\partial B(z,r)}\Sigma (x).
\eeq
Then there exists $\eps_{0}>0$ such that, for every 
$0<\eps<\eps_{0}$, problem $(P_{\eps})$ admits a nontrivial  
solution $u_{\eps}\in\H$, $u_{\eps}\geq0$, such  
that the following facts hold:
\vskip3pt
\noindent(i) $u_{\eps}$ admits exactly one global maximum
point $x_{\eps}\in B(z,r)$;
\vskip3pt
\noindent(ii) $\dys \lim_{\eps\to 0} \Sigma(x_{\eps})=\Sigma_0;$
\vskip3pt
\noindent(iii) there exist $\mu_1,\mu_2>0$ such that, 
for every $x\in\R^N$,
$$
u_\eps(x)\leq \mu_1e^{-\mu_2\textstyle\frac{|x-x_\eps|}{\eps}}.
$$
\ete

\bos
In this abstract result
\[
\Sigma_0=\Sigma(z)=I_z(Q_z),
\]
where $Q_z$ is the unique positive least energy critical point
(see \cite{stzo,seta,sezo}) of $I_z$ defined in \eqref{Iy}

 $V_{0}$ and $s_{0}$ are just the value 
of the functions $s(x)$ and $V(x)$ on $z$, and they are not in general 
related with the minimum values of $s(x)$ and $V(x)$, because
it is actually the minimum point of the function $\Sigma$ that
plays the fundamental role. In Theorem \ref{mainvs} we have seen
that, in the particular case in which $V$ and $s$ attain their
minimum in the same point, then this  point will be 
a minimum of $\Sigma$. But, in general, 
this could not be the case, and
still we may find a minimum point of $\Sigma$.
\eos

\bos
In order to find a sufficient condition in terms of 
an explicit concentration function, instead of $\Sigma$, 
it is usually crucial to find a change of variable, from
the frozen problem to the problem with with all the 
constants equal to one, that is in this case from
\eqref{limity} to 
$$
\begin{cases}
\dys -\Delta u+u^{2}=\frac{u^{3}}{1+u^{2}} &\text{in } \R^{N},
\\
\dys \,u\,\in H^{1}(\R^{N}).
\end{cases}
$$
Unfortunately, this procedure is feasible in this context. 
Since doing that we move from a problem which admits solutions 
to one which has not any nontrivial solution.
This implies that we cannot express a solution as a member 
of a two-parameters family generating by a fundamental solution
as in the most studied case \cite{wangzeng,amfema,bova}
\eos

With respect to the topic of locating the possible concentration
points, let us first introduce the concentration set. 

\begin{definition}\label{defe}
The concentration set $\mathcal{E}$ for 
problem~\eqref{Eq:eps}, is defined by
\begin{align*}
\mathcal{E}=\Big\{& z\in\R^N \hbox{ such that there exists 
a sequence of solutions $\{u_\eps\}\in\H$ of \eqref{Eq:eps} with}\\
& \text{$u_{\eps}(z+\eps x)\to 0$ as $|x|\to\infty$ uniformly w.r.t.
 $\eps$ and ${\eps}^{-N}J_{\eps}(u_{\eps})\to\Sigma(z)$
as $\eps\to 0$} \Big\},
\end{align*}
where $\Sigma$ is defined in \eqref{defsigmagen}.
\end{definition}

We will prove the following result concerning necessary 
conditions for the concentration to occur.
\bte\label{mainnece}
Assume \eqref{s0}, \eqref{V0} and that $V,\,s\in C^{1}(\R^{N})$ 
such that there exist  $\beta>0$, $\gamma\geq 0$, satisfying
\beq\label{growths}
|\nabla V|\leq \beta e^{\gamma |x|} \quad \text{and }\quad
|\nabla s|\leq \beta e^{\gamma |x|}\quad \forall \,x\in \R^{N}.
\eeq
Then $\Sigma $ is of class $C^{1}(\Omega)$ and if $z\in {\mathcal E}$
the following facts hold:
\\
(i)\;  $\nabla V(z)$ and $\nabla s(z)$ are linearly dependent and point
in opposite directions.
\\
(ii)\; Either $\partial_{j} V(z)=\partial_{j} s(z)=0$ for every $j=1,\dots, N$
or  there exists at least a $j_{0}\in \{1,\dots, N\}$ such that
$\partial_{j_{0}} V(z),\,\partial_{j_{0}} s(z)\neq 0$ with
still $\partial_{j} \Sigma(z)=0$ for every $j=1,\dots, N$.
\ete
\bos\label{quasifine}
We can say a little bit more in the last conclusion of the above Theorem.
Indeed, as a consequence of conclusion {\it (i)}, every nontrivial
partial derivative of $V$ and $s$ satisfies a precise identity (see for more
details Remark \ref{fine}).
\eos

\section{Proofs of Theorems \ref{mainvs} and \ref{mainsigma}}\label{proof}

In this section we will prove Theorems \ref{mainvs} 
and \ref{mainsigma} using the well known penalization 
procedure introduced in \cite{defe,defe2}. Let us deal first 
with the proof of Theorem \ref{mainvs}.
\\
Recalling that the derivative of $F(x,t)$ with respect 
to $t$ is given by $f(x,t)=t^{3}/(1+s(x)t^{2})$, fix 
$0<\nu<1/2$ and define the function 
$\ov{f}(x,t):\R^{N}\times \R\to \R$ by
\beq\label{fbar}
\ov{f}(x,t)=
\begin{cases}
\min\left\{f(x,t),\nu \mu t\right\} & \text{ for } t\geq 0,
\\
0 &\text{ for } t<0.
\end{cases}
\eeq
Let $0<r'<r$ such that 
\beq\label{vsnew}
\left\{\begin{array}{c}
\dyle V_{0}=V(z)=\min_{B(z,r')}V (x)\leq\min_{\partial B(z,r')}V (x), 
\quad
s_{0}=s(z)=\min_{B(z,r')}s (x)<\min_{\partial B(z,r')}s (x),\\ 
\hbox{or}\\
\dyle V_{0}=V(z)=\min_{B(z,r')}V (x)<\min_{\partial B(z,r')}V (x), 
\quad
s_{0}=s(z)=\min_{B(z,r')}s (x)\leq\min_{\partial B(z,r')}s (x).\\
\end{array}\right.
\eeq

Note that the existence of $r'$ is a consequence of the
continuity of the functions $s,V$ in $\R^N$ and of \eqref{minivs}.
Indeed, let us prove the existence of $r'$ in the case in which the
first assumption in \eqref{minivs} is satisfied.
Arguing by contradiction, it follows that for any $\rho<r$ it holds
\[
\inf_{B(z,\rho)}V (x)>\min_{\partial B(z,\rho)}V (x)
\hbox{ or }
\inf_{B(z,\rho)}s (x)\geq\min_{\partial B(z,\rho)}s (x).
\]
Since, for every continuous function it holds
\[
\inf_{B(z,\rho)}k (x)=\min_{\ov{B(z,\rho)}}k(x)
\leq \min_{\partial B(z,\rho)}k (x)
\]
the first inequality for $V$ cannot be true and
we can reduce to the case
\[
\inf_{B(z,\rho)}s (x)=\min_{\partial B(z,\rho)}s (x)
\qquad \forall \rho<r.
\]
Now, let $\{\rho_n\}$ be an increasing sequence such that
$\rho_n \to r$. Then we can write
\[
\inf_{B(z,\rho_n)}s (x)=\min_{\ov{B(z,\rho_n)}}s(x)
= \min_{\partial B(z,\rho_n)}s (x)=s(p_n)
\]
with $p_n\in\partial B(z,\rho_n)$.

As $\partial B(z,\rho_n)\subset\overline{B(z,r)}$
and $\rho_n\to r$, it results, up to a subsequence, 
$p_n\to p\in \partial B(z,r)$ and, passing to the limit,
\[
\min_{\ov{B(z,r)}}s(x)
=\lim_{n\to+\infty}\min_{\ov{B(z,\rho_n)}}s(x)
=\lim_{n\to+\infty}s(p_n)=s(p)
\geq \min_{\partial B(z,r)}s(x)
\]
and this contradicts \eqref{minivs}. This proves the existence
of $r'<r$ such that \eqref{vsnew} holds.

Let $\chi\in C^\infty(\R^N)$ be such that
\beq\label{defchi}
\dyle \chi(x)= 1 \quad \forall x\in B(z,r'),  \quad
\dyle \chi(x)= 0 \quad \forall x\in \R^N \setminus \ov{B(z,r)},
\eeq
and set
\[
g(x,t)=\chi(x) f(x,t)+(1-\chi(x))\ov{f}(x,t),
\]
for a.e. $x\in\R^N$ and any $t\in\R$.
Having defined $ G(x,t)=\int_{0}^{t}g(x,\xi)d\xi$, 
in the light of the above definition, the following
result follows.
 
\ble\label{NQpernoi}
Assume conditions \eqref{minivs}. Then, the following conditions
hold for every $t$ in $\R$ and for almost every $x$ in $\R^{N}$
\beq\label{gzero}
\lim_{t\to0}\frac{g(x,t)}{t}=0, \quad\text{uniformly in } x\in \R^{N}
\eeq
\begin{align}\label{NQ1}
g(x,t) t-2G(x,t)&\geq 0 ,\;\; \lim\limits_{|t|\rightarrow \infty}  
g(x,t)t - 2G(x,t) = +\infty, 
 \quad\;\; \forall x\in B(z,r),\,
\\
\label{NQ2}
0 &\leq 2G(x,t) \leq  g(x,t)t\leq \nu \mu t^{2} 
\qquad\qquad\quad \forall x\not\in\overline{B(z,r)}
\end{align}
\ele
\bdim
For $t$ going to zero, $\ov{f}(x,t)=f(x,t)$ and from \eqref{s0} it results
$$
\frac{g(x,t)}{t}=\frac{f(x,t)}t=\frac{t^{2}}{1+s(x)t^{2}}\leq t^{2},
$$
implying \eqref{gzero}.
\\
In order to show \eqref{NQ1} and \eqref{NQ2}, let us first show
that \eqref{nonq} holds because
$$
\frac12 f(x,t)t-F(x,t)= \frac1{2s^{2}(x)}\left[\ln(1+s(x)t^{2})
-\frac{s(x)t^{2}}{1+s(x)t^{2}}\right].
$$
Then \eqref{NQ1} easily follows studying the function
$h(t)=\ln(1+t)-t/(1+t)$.

Now, if $x\in B(z,r)$ then 
\[
g(x,t)t-2G(x,t)\geq \chi(x)\left[ f(x,t)t-2F(x,t) \right]
\]
so that, for every $x\in B(z,r)$, \eqref{nonq} implies
\eqref{NQ1} being $\chi(x)>0$. 

For $x\not\in \ov{B(z,r)}$, $g(x,t)=\ov{f}(x,t)$ and
$$
G(x,t)=\begin{cases}
F(x,t) &\text{ if } f(x,t)\leq \nu \mu t
\smallskip\\
\dys \frac{\nu \mu}2 t^{2} & \text{ if } f(x,t)> \nu \mu t
\end{cases}
$$
and \eqref{NQ2} easily follows.
\edim
The following easy technical lemma will be useful in the sequel.
\ble\label{decre}
The following facts hold:\\
i) For any $2\leq q\leq 4$, there exists $C=C(q)$ such that
 $t^2-\ln(1+t^2)\leq C|t|^{q}$,  in $\R$. \\
ii) $t^{2}/(1+st) \leq C|t|$, for all $t$ in $\R$ 
and  for every $s\in \R^{+}$.\\
iii) For every $L\geq 0$, the real function 
$h\in C(\R^{+},\R^{+})$ defined
$\dys h(s)=\frac{L}{s}-\frac1{s^{2}}\ln(1+Ls)$
is monotone decreasing.
\ele
\bdim
The proof can be shown by direct calculations.
\edim

We will study the penalized
 functional $J_{\eps}:\H(\R^{N})\to \R$
defined by
\beq\label{jeps}
J_{\eps}(u)= \frac{1}{2}\| u\|_{\eps,V}^{2}
-\intr G(x,u(x))dx
\eeq
whose critical points are solutions of the problem
\beq\label{pbareps}
\begin{cases}
\dys -\eps^{2}\Delta u+V(x)u=g(x,u(x)) &\text{in } \R^{N},
\\
\dys \,u\,\in {\H}(\R^{N}).
\end{cases}
\eeq
In the study of asymptotically linear problems the usual Palais-Smale
condition is substituted by the following Cerami condition introduced in \cite{cerami}.
\begin{definition}\label{defcer}
Let $E$ be a Banach space.
A sequence $\{u_{n}\} \subset E$ is said to be a Cerami sequence 
for a functional $I$, $(Ce)_c$ for short, if
\beq\label{Iconv} 
I(u_n) \to c, \qquad
(1+\left\|u_n \right\|_{E})  \left\| I'(u_n) \right\|_{E^*} \to 0 .
\eeq
Moreover, a functional  $I \in C^1(E,\R)$ is said to satisfy 
the Cerami condition $(Ce)_c$ if any Cerami sequence  
possesses a convergent subsequence.
\end{definition}
In the next lemma we prove that $\jbe$ satisfies the Cerami condition.

\ble\label{ceramicomp}
Assume conditions  \eqref{s0}, \eqref{V0} and \eqref{minivs}. 
Then, for every $\eps>0$ fixed, every Cerami sequence for 
$\jbe$ admits a convergent subsequence.
\ele

\bdim Let us take $\{u_{n}\}$ a Cerami sequence and let us 
first prove that $\{u_{n}\}$ is bounded by contradiction.
Assume then, up to a subsequence, 
\beq\label{cer}
\|u_{n}\|_{\eps,V}\to \infty, \quad \jbe(u_{n})\to c_{\eps},\quad 
\|\jbe'(u_{n})\|_{(\H(\R^{N}))^{*}}\|u_{n}\|_{\eps,V}<\frac1n.
\eeq
Arguing as in Lemma 3.30 of \cite{mamopesat} it is possible to obtain the 
following inequalities for every $t$
\begin{align}
\label{nqg}
\intr \frac{1}{2}g(x,u_n(x)) u_n(x) - G(x,u_n(x)) 
&\leq \jbe(u_n) +\opic(1) \leq c_{\eps}+\opic(1),
\\ 
\label{gpos}
\intr G(x,t_{n}u_{n}(x))
&\geq 
\frac{t^{2}}{2}\|u_{n}\|_{\eps,V}^{2}-\eps^{N}c_{\eps}+ \opic(1),
\\
\label{stimeJ}
J_\eps(tu_{n}) &\leq J_\eps(u_{n}) +\opic(1),
\end{align}
with $\opic(1)\to0$ as $n\to\infty$.

Let us define the  function
$\phi_{n}(x)=u_{n}(z+\eps x)$
and notice that $\phi_{n}$ belongs to the Hilbert space
\beq\label{hez}
\varmathbb{H}_{\eps,z}\,:\,=\left\{v\in H^{1}(\R^{N})\,:\,
\intr V(z+\eps x)v^{2}(x)dx<+\infty \right\}
\eeq
with norm 
$$
\|v\|^2_{\eps,z}:=\|\nabla v\|_{2}^{2}+\intr V(z+\eps x)v^{2}(x)dx.
$$
This class of spaces has been used in \cite{jeta2} and
in the rest of this proof (and in the proof of Lemma \ref{tresei})
we will adopt some of their arguments.
First, notice that proving the result is equivalent 
to show that  there exists a positive constant $C$ 
possibly depending on
$\eps$, satisfying $\|\phi_{n}\|_{\eps,z}^{2}\leq C$
for every $\eps$ fixed, as it results
\beq\label{normafi}
\|\phi_{n}\|_{\eps,z}^{2}=
\intr \left[\eps^{2}|\nabla u_{n}(z+\eps x)|^{2}+
V(z+\eps x)u_{n}^{2}(z+\eps x)\right]dx
=\frac{\|u_{n}\|_{\eps,V}^{2}}{\eps^{N}}.
\eeq
Let us study  \eqref{gpos} in terms of 
$\phi_{n}$. Since 
$$
\intr G(x,tu_{n}(x))dx=\eps^{N}\intr G(z+\eps \xi, tu_{n}(z+\eps \xi))
d\xi=\eps^{N}\intr G(z+\eps \xi, t\phi_{n}(\xi))d\xi,
$$
the sequence $\phi_{n}$ satisfies the following inequality
$$
\intr G(z+\eps \xi,t\phi_{n}(\xi))d\xi\geq 
\frac{t^{2}}{2}\|\phi_{n}\|_{\eps,z}^{2}-c_{\eps}+ \opic(1),
$$
and choosing $t=t_{n}=2\sqrt{c_{\eps}}/\|\phi_{n}\|_{\eps,z}$,
it follows
\[
\intr G(z+\eps \xi,t_{n}\phi_{n}(\xi))d\xi\geq c_{\eps}+ \opic(1).
\]
Now, we can argue by contradiction, supposing 
that, up to a subsequence, $\|\phi_{n}\|_{\eps,z}\to+\infty$ and
defining the sequence $\psi_{n}=t_{n}\phi_{n}$, which
verifies $\|\psi_{n}\|_{\eps,z}=2\sqrt{c_{\eps}}$, 
to obtain $\psi\in H^{1}$ such that $\psi_{n}$ converges to
$\psi$ weakly in $H^{1}$, strongly in  $L^p_{{\rm loc}}(\RN)$
$\forall\, p\,\in [1,2^*)$, and almost everywhere.
We claim that
\begin{equation}\label{pallafissa}
\limsup_{n\to 0} \sup_{y\in \R^{N}}\int_{B(y,1)}
|\chi(z+\eps x)\psi_n(x)|^2 dx>0,
\end{equation}
where $\chi$ is introduced in \eqref{defchi}. By contradiction, 
if \eqref{pallafissa} were false, it would result
$$
\lim_{n\to 0}\sup_{y\in \RN}\int_{B(y,1)}|\chi(z+\eps x)
\psi_n(x)|^2dx= 0 ,
$$
then, using the argument of Lemma I.1 in \cite{plions}, we deduce 
that, for $2<p<2^*$, the sequence $\chi(z+\eps x)\psi_n$ converges 
to $0$ strongly in $ L^{p}(\RN)$.
Let us fix $q\in(2+2/N,\ov{q})$, with $\ov{q}=\min\{4, N/(N-2)+1\}$
and apply conclusion {\it i)} in Lemma \ref{decre},
H\"older inequality and \eqref{defchi} 
to obtain for every $L\geq 1$
\begin{align}\label{stimaG}
\intr|\chi(z+\eps x)F(z+\eps x,L\psi_{n}(x))|dx &\leq C
L^{q}\intr\chi(z+\eps x)|\psi_{n}(x)|^{q}dx
\\
\nonumber &\leq 
C L^{q}
\left[\intr\chi^p(z+\eps x)|\psi_{n}|^{p}\right]^{1/p}
\|\psi_{n}\|_{p'(q-1)}^{(q-1)},
\end{align}
where $p'=p/(p-1)$. Since $q\in (2+2/N,\ov{q})$, 
$p'(q-1)\in (2,2^{*})$, so that the last
integral is uniformly bounded, implying that
$\chi(z+\eps x)F(z+\eps x,L\psi_{n}(x))$ converges to zero
in $L^{1}(\R^{N})$. Moreover, for every $L\geq 1$, 
from \eqref{fbar} and \eqref{V0} it follows
$$
\intr (1-\chi(z+\eps x))\overline{F}(z+\eps x,L\psi_{n})
\leq  L^{2} \frac{\nu\mu}2\intr |\psi_{n}|^{2}
\leq \frac{\nu L^{2}}2
\|\psi_{n}\|_{\eps,z}^{2}\leq c_\eps L^{2}.
$$
Therefore for every $L\geq 1$
\beq\label{disL}
\tilde{J}_{\eps}\left(L \psi_{n}\right)
\geq c_\eps L^{2} +\opic(1),
\eeq
where $\tilde{J}_{\eps}:\varmathbb{H}_{\eps,z}\to \R$ 
is defined by
\beq\label{defjtilde}
\tilde{J}_{\eps}(v):=\frac12 \|v\|_{\eps,z}^{2}-\intr G(z+\eps x,v(x))dx.
\eeq
On the other hand, from \eqref{stimeJ} we deduce that
$$
\tilde{J}_{\eps}(L\psi_{n})
=\tilde{J}_{\eps}\left(L t_n \phi_{n}\right)
=\frac1{\eps^{N}} J_{\eps}(Lt_{n}u_{n}) \leq
\frac1{\eps^{N}}\left[J_{\eps}(u_{n})+\opic (1)\right]
\leq \frac1{\eps^{N}}\left[ c_{\eps}+\opic(1)\right].
$$
This together with \eqref{disL} produce a contradiction,
yielding \eqref{pallafissa}.
\\
As in \cite{jeta2}, this implies the existence of a number
$\gamma>0$, and of a sequence $\{y_{n}\}$ with 
$B(y_{n},1)\cap$supp$\chi(z+\eps \cdot)\neq\emptyset$ 
and such that
\beq\label{pallafissa2}
\lim_{n\to +\infty}\int_{B(y_{n},1)}
|\chi(z+\eps x)\psi_n(x)|^2 dx >0.
\eeq
Since $B(y_n,1)\cap\hbox{supp}\chi(z+\eps x)\neq\emptyset$,
\eqref{defchi} implies that
there exists a $\eta$ satisfying
$|y_n-\eta|<1 $ and $|z+\eps\eta-z|<r,$
so that $|\eps y_n|\leq \eps|y_n-\eta|+\eps|\eta|<\eps+r$
and  we can find  $x_0$ such that 
\beq\label{xzerofissa}
\eps y_{n}\to x_{0}\in \ov{B(0,r+\eps)}.	
\eeq
Let us now define the functions
$$
\ov{\psi}_{n}(x)=\psi_{n}\left(y_{n}+x\right),\qquad 
\ov{\chi}_{n}(x)=\chi\left(z+\eps\left(y_{n}+x\right)\right),
$$
and observe that, as $\ov{\psi}_{n}$ is uniformly bounded in 
$H^{1}(\R^{N})$, it exists $\ov{\psi}\in H^{1}(\R^{N})$ such that
$\ov{\psi}_n$ converges to $\ov{\psi}$ weakly in $H^{1}(\R^{N})$, 
 almost everywhere and strongly in $L^{2}(B(0,1))$.
Moreover, from \eqref{xzerofissa} we deduce that
$\ov{\psi}_{n}(x)\ov{\chi}_{n}(x)\to \ov{\psi}(x)\chi(z+x_{0}+\eps x)$
almost everywhere, then \eqref{pallafissa2} yields 
\begin{align*}
0 <  \lim_{n\to \infty}\int_{B(y_{n},1)}|\chi(z+\eps \xi)\psi_n(\xi)|^2dx & =\lim_{n\to \infty}\int_{B(0,1)}{\Big |}
\psi_{n}\big(y_{n}+x\big)
\chi\big(z+\eps\big(y_{n}+x\big)\big)
{\Big |}^{2}dx \\
=\lim_{n\to \infty}\int_{B(0,1)}|\ov{\psi}_{n}(x)
\ov{\chi}_{n}(x)|^{2}dx
&=\int_{B(0,1)} \chi^{2}(z+x_{0}+\eps x) |\ov{\psi}(x)|^{2}dx
\end{align*}
which implies that there exists
an open set $A\subset B(0,1)$ such that for every $x\in A$ it holds
$|\ov{\psi}(x)|>0$ and $\chi(z+x_{0}+\eps x)>0$. Moreover, it results
$$
0<|\ov{\psi}(x)|=\lim_{n\to \infty}\Big|
\psi_{n}\left(y_{n}+x\right)\Big|
=2\sqrt{c_{\eps}}\lim_{n\to \infty}
\frac{\Big|\phi_{n}\left(y_{n}+ x\right)\Big|}
{\|\phi_{n}\|_{\eps,z}}
$$
then, for every $x\in A$, we have that
$\big|\phi_{n}(y_{n}+x)\big|\to +\infty$ and
as $z+\eps(y_{n}+x) \to z+x_0 +\eps x$,
with $\chi(z+x_{0}+\eps x)>0$, \eqref{defchi}
yields the existence of $n_{0}$ 
such that, for $n\geq n_{0}$, 
$z+\eps( y_{n}+x)\in B(z,r)$. Then,  
\eqref{NQ1} and \eqref{NQ2} give
$$
\begin{array}{c}
\dys \lim_{n\to \infty} \int_{{\R^{N}}} 
\frac{1}{2}g(z+\eps(y_{n}+x),
\phi_{n}(y_{n}+x))\phi_{n}(y_{n}+x)
\!-\!G(z+\eps(y_{n}+x), \phi_{n}(y_{n}+x))
\geq 
\medskip\\
\dys \lim_{n\to \infty} \int_{A} 
\frac{1}{2}g(z+\eps (y_{n}+x),\phi_{n}(y_{n}+x))
\phi_{n}(y_{n}+x)
-G(z+\eps(y_{n}+x),\phi_{n}(y_{n}+x))=+\infty
\end{array}
$$
But, on the other hand, performing the change of variable 
$z+\eps(y_{n}+x)=\xi$, and using 
 \eqref{nqg},  one derives the desired contradiction.
\edim
\bos
In the proof of the above lemma we use 
assumption \eqref{minivs} only to define the penalization
with the function $g(x,t)$ satisfying condition
\eqref{NQ1}, \eqref{NQ2}. 
\eos

\ble\label{mountaineps}
Assume \eqref{s0}, \eqref{V0} and \eqref{minivs}. Then
there exists $\eps_{0}>0$ such that, for every $\eps\in(0,\eps_{0}) $ 
$J_{\eps}$ has a nontrivial critical point
$u_{\eps}$ satisfying
\beq\label{leveleps}
 J_{\eps}(u_{\eps})\leq \eps^{N}(\Sigma(z)+\opic(1)),
\eeq
where  $\opic(1)\to 0$ as $\eps\to 0$. 
\ele
\bdim
We will obtain the existence of $u_{\eps}$ by applying the 
variant of the Mountain Pass Lemma with the Cerami condition
(see \cite{cerami, babefo}) 
to the functional ${J}_{\eps}$. Let us first notice that  
condition \eqref{gzero} immediately implies that $v_{0}=0$ is a strict 
local minimum. In order to show the existence of $v^{*}$ 
such that ${J}_{\eps}(v^{*})<0$, 
let us observe that, arguing as in Lemma 2.1 \cite{jeta1}, 
we find $w^{*}$ such that $I_{z}(w^{*})<0$,
for $I_{z}(v)$ defined in \eqref{Iy}.
Let us choose $v^{*}_{\eps}(x)=\eta(x)w^{*}\big(\frac{x-z}{\eps}\big)$, with
$\eta(x)$  a smooth function compactly supported in $\R^{N}$ and 
$\eta(x)\equiv 1$ in $B(z,r)$.
Computing $\jbe(v^{*}_{\eps})$ gives
\begin{align*}
\jbe(v^{*}_{\eps})&=\!\!\!\frac{\eps^{2}}{2}\!\intr \!\!\!
\left\{
|\nabla \eta|^{2}|w^{*}\big(\frac{x-z}{\eps}\big)|^{2}+
\frac2{\eps}\eta
w^{*}\big(\frac{x-z}{\eps}\big)\nabla \eta\nabla w^{*}\big(\frac{x-z}{\eps}\big)
+\frac1{\eps^{2}}\eta^{2} |\nabla w^{*}\big(\frac{x-z}{\eps}\big)|^{2}
\right\}\!\!dx
\\
&+\frac{1}2\intr V(x) \eta^{2}|w^{*}\big(\frac{x-z}{\eps}\big)|^{2}dx
-\intr G\left(x,\eta w^{*}\big(\frac{x-z}{\eps}\big)\right)dx,
\end{align*}
and performing the change of variable $y=(x-z)/\eps$, 
using the properties of the function $\eta$, one gets
$$
\jbe(v^{*}_{\eps}) =
\dys \eps^{N}\intr \left[\frac12\eta^{2}(z+\eps y) \big(|\nabla w^{*}|^{2}
+V(z+\eps y) |w^{*}|^{2}\big)-G\left(z+\eps y,\eta(z+\eps y) w^{*}\right)\right]+\text{ o}\big(\eps^{N}\big),
$$
where o$(\eps^{N})/\eps^{N}\to 0$ as $\eps$ goes to zero.
Since the above integral uniformly converges, 
as $\eps$ goes to zero, to $I_{z}(w^{*})<0$, 
there exists $\eps_{0}>0$  such that for $\eps\in (0,\eps_{0})$,
$\jbe(v^{*}_{\eps})<0$, giving the desired conclusion.
\\
The geometric behavior  just observed yields  the construction of
a Cerami sequence $\{u_{n}\}$ of $\jbe$  at the Mountain Pass level
$c_{\eps}$, defined by
\beq\label{defgamma}
c_{\eps}=\inf_{\gamma\in\Gamma_{\eps}}
\max_{[0,1]}J_{\eps}(\gamma(t)),\quad 
\Gamma_{\eps}=\big{\{}\gamma\in C([0,1],\H(\R^{N})),\;
\gamma(0)=0,\, J_{\eps}(\gamma(1))<0 {\big \}}.
\eeq
Then, Lemma \ref{ceramicomp} allows to pass to the limit
and obtain a critical point $u_{\eps}$ with 
$J_{\eps}(u_{\eps})=c_{\eps}$.
\\
In order to show \eqref{leveleps}, we will argue as in 
Proposition 6.1 in \cite{jeta2}. From Lemma 2.1 in \cite{jeta1}
we deduce the existence of a path $\gamma\in C([0,1],H^{1}(\R^{N}))$
such that
\beq\label{livello}
\gamma(0)=0,\quad  I_{z}(\gamma(1))<0,\quad 
I_{z}(\gamma(t))\leq \Sigma(z),\quad
\max_{t\in [0,1]}I_{z}(\gamma(t))=I_{z}(Q_{z})=\Sigma(z),
\eeq
where $Q_{z}$ is the unique positive solution of \eqref{limity}
with $y=z$ (see Section \ref{setting}) and 
$\Sigma(z)$ is defined in \eqref{defsigmagen}.

Let us consider a function $\eta\in C^{\infty}_{0}(\R^{N})$ 
such that $\eta(0)=1$ and $0\leq\eta(x)\leq 1$. 
First, notice that there exists $R_{0}$ such that 
\beq\label{nega}
I_{z}(\eta(z/R)\gamma(1))<0, \quad \forall \, R\geq R_{0}
\quad \hbox{and}\quad
\lim_{R\to +\infty}\max_{t\in [0,1]}
I_{z}(\eta(z/R)\gamma(t))=  \Sigma(z)
\eeq

Then we define the path
$$
\gamma_{R,\eps}(t)(y)=\eta\left(\frac{y}{R}
\right)\gamma(t)\left(\frac{y-z}{\eps}\right)
$$
so that $\gamma_{R,\eps}(t):[0,1]\to H^{1}(\R^{N})$.
Since $J_{\eps}\left(\gamma_{R,\eps}(t)\right)/\eps^{N}$ 
converges to $ I_{z}\big(\eta(z/R)\gamma(t)\big)$ 
as $\eps$ goes to zero, uniformly with respect to 
$t\in [0,1]$,
\eqref{nega} implies that $\gamma_{R,\eps}\in \Gamma_{\eps}$.
Moreover, from \eqref{nega} it follows
$$
c_{\eps}\leq \max_{[0,1]}\frac{J_{\eps}
\left(\gamma_{R,\eps}(t)
\right)}{\eps^{N}}=\max_{[0,1]}
I_{z}\left(\eta(z/R)\gamma(t)\right)+\opic(1)
\leq \Sigma(z)+\opic(1),
$$
implying \eqref{leveleps}.
\edim

\ble\label{tresei}
There exists a positive constant $L$ and $\eps_{0}>0$
such that for every $\eps\in (0,\eps_{0})$
\beq\label{stimato}
\|u_{\eps}\|_{\eps,V}^{2}\leq L\eps^{N}.
\eeq
\ele

\bdim We will follow the argument of Lemma \ref{ceramicomp} paying
attention to the fact that now $\eps$ is not fixed.
Arguing as in Lemma 3.30 of \cite{mamopesat} it is possible to obtain 
the following inequalities  for every $t$
\begin{align}
\label{stime}
J_{\eps}(tu_{\eps}) \leq J_{\eps}(u_{\eps}),
\quad
\intr G(x,tu_{\eps}(x))\geq \frac{t^{2}}{2}
\|u_{\eps}\|_{\eps}^{2}-\eps^{N}(\Sigma(z)+ \opic(1)).
\end{align}
Introducing the function $\phi_{\eps}(x)=u_{\eps}(z+\eps x)$
belonging to $\varmathbb{H}_{\eps,z}$, defined in \eqref{hez},
notice that, \eqref{normafi} tells us that proving \eqref{stimato} 
is equivalent to show that there exists $\eps_{0}>0$ such that  
$\|\phi_{\eps}\|_{\eps,z}^{2}\leq L$ for every $\eps\in(0,\eps_{0})$.

As in Lemma \ref{ceramicomp} we obtain that the sequence $\phi_\eps$
satisfies the following inequality
\beq\label{stimaphi}
\intr G(z+\eps \xi,t_{\eps}\phi_{\eps}(\xi))d\xi\geq \Sigma(z)+ \opic(1),
\eeq
with $t_\eps=2\sqrt{\Sigma(z)}/\|\phi_\eps\|_{\eps,z}$.
Arguing again by contradiction and supposing 
that, up to a subsequence, $\|\phi_{\eps}\|_{\eps,z}\to+\infty$,
we set  $\psi_{\eps}=t_{\eps}\phi_{\eps}$, which
verifies $\|\psi_{\eps}\|_{\eps,z}=2\sqrt{\Sigma(z)}$. Then
$\psi_{\eps}$ converges weakly in $H^{1}$ strongly in 
$L^p_{{\rm loc}}(\RN)$, $\forall\, p\,\in [1,2^*)$, almost everywhere
and, as in Lemma \ref{ceramicomp}, it satisfies
\begin{equation}\label{palla}
\limsup_{\eps\to 0} \sup_{y\in \R^{N}}\int_{B(y,1)}
|\chi(z+\eps x)\psi_\eps(x)|^2 dx>0.
\end{equation}

This implies the existence of a sequence $\{y_{\eps}\}$ with
$B(y_\eps,1)\cap\hbox{supp}\chi(z+\eps\cdot)\neq\emptyset$, and
such that 
\beq\label{palla2}
\lim_{\eps\to0}\int_{B(y_{\eps},1)}
|\chi_{\eps}(z+\eps x)\psi_\eps(x)|^2 dx >0.
\eeq
Moreover, from \eqref{defchi}, it follows that 
$\eps y_\eps\in\{\eta:|\eta|<\eps+r\}$ so that
\beq\label{xzero}
\eps y_{\eps}\to x_{0}\in \ov{B(0,r)}.
\eeq
In this case, the functions
$$
\ov{\psi}_{\eps}(x)=\psi_{\eps}\left(y_{\eps}+x\right),\qquad 
\ov{\chi}_{\eps}(x)=\chi\left(z+\eps\left(y_{\eps}+x\right)\right)
$$
are such that $\ov{\psi}_{\eps}$ converges to $\ov{\psi}$ 
weakly in $H^{1}(\R^{N})$, almost everywhere and strongly 
in $L^{2}(B(0,1))$. 
Moreover, from \eqref{xzero} we deduce that
$\ov{\psi}_{\eps}(x)\ov{\chi}_{\eps}(x)\to \ov{\psi}(x)\chi(z+x_{0})$
almost everywhere, then \eqref{palla2} yields 
$$
0<\lim_{\eps\to 0}\int_{B(0,1)}|\ov{\psi}_{\eps}(x)\ov{\chi}_{\eps}(x)|^{2}dx =\chi^{2}(z+x_{0})\int_{B(0,1)}|\ov{\psi}(x)|^{2}dx.
$$
Then $\chi(z+x_{0})>0$ implying that $x_{0}\in B(0,r)$
and there exists an open set $A\subset B(0,1)$ such that for every 
$x\in A$, $|\ov{\psi}(x)|>0$. Moreover, it results
$$
0<|\ov{\psi}(x)|=\lim_{\eps\to 0}\Big|\psi_{\eps}\left(y_{\eps}+x\right)\Big|=4\Sigma(z)\lim_{\eps\to 0}\frac{\Big|\phi_{\eps}\left(y_{\eps}+x\right)\Big|}
{\|\phi_{\eps}\|_{\eps,z}}
$$
so that, for every $x\in A\subset B(0,1)$, 
$\big|\phi_{\eps}(y_{\eps}+x)\big|\to +\infty$.
Moreover, as $z+\eps(y_{\eps}+x)\to 
z+x_{0}\in B(z,r)$ 
we can deduce that, for $\eps$ sufficiently small, 
$z+\eps (y_{\eps}+x)\in B(z,r)$, 
so that \eqref{NQ1} and \eqref{NQ2} give
\begin{align*}
\lim_{\eps\to 0}\intr \dfrac{1}{2}g(z+\eps(y_{\eps}+x),
\phi_{\eps}(y_{\eps}+x))\phi_{\eps}(y_{\eps}+x)
-G(z+\eps(y_{\eps}+x),\phi_{\eps}(y_{\eps}+x))dx
&\geq
\\
\lim_{\eps\to 0}\int_{A} 
\dfrac{1}{2}g(z+\eps(y_{\eps}+x),
\phi_{\eps}(y_{\eps}+x))\phi_{\eps}(y_{\eps}+x)
-G(z+\eps(y_{\eps}+x),\phi_{\eps}(y_{\eps}+x))
 dx &=+\infty.
\end{align*}
But, on the other hand, it results
$$
\intr \left[\dfrac{1}{2}g(z+\eps(y_{\eps}+x),
\phi_{\eps}(y_{\eps}+x))\phi_{\eps}(y_{\eps}+x)
-G(z+\eps\xi_{\eps},\phi_{\eps}(y_{\eps}+x))
\right]dx= \tilde{J}_{\eps}(\phi_{\eps})
$$
which is uniformly bounded because of  \eqref{leveleps}.
\edim

\begin{proposition}\label{maximo}
Assume  \eqref{s0}, \eqref{V0} and \eqref{minivs}. 
Then for every $\delta>0$
there exists $\eps_\delta>0$ such that
\begin{equation}\label{limepsfin}
\sup_{0<\eps<\eps_\delta}\,\sup_{x\in\R^{N}\setminus B(z,r)}u_{\eps}(x)<\delta.
\end{equation}
\end{proposition}
\bdim
Let us first prove that
\begin{equation}\label{limeps}
\lim_{\eps\to 0}\,\sup_{x\in\partial B(z,r)}u_{\eps}(x)=0.
\end{equation}
We proceed by contradiction, assuming that there exist
a sequence $\{\eps_{n}\}$ con\-ver\-ging to $0$ and a sequence 
$\{x_{n}\}\subset\pa B(z,r)$ such that, for some positive constant $\beta$,
\beq\label{contro}
u_{\eps_n}(x_{n}) \geq \beta\qquad\text{for all $n\geq 1$.}
\eeq
Since $\pa B(z,r)$ is a compact set, we can assume that there exists a
subsequence of $\{x_{n}\}$, still denoted by $\{x_{n}\}$, which converges to a
point $x_{0}\in\pa B(z,r)$. Consider the scaling of $u_{\eps_{n}}$ 
centered at $x_n$, that is
$$
\phi_{n}(x)=u_{\eps_{n}}(x_{n}+\eps_{n}x),
$$
which  solves the equation
\beq\label{succ}
-\Delta\phi_{n}+V(x_{n}+\eps_{n}x)\phi_{n}=g(x_{n}+\eps_{n}x,\phi_{n}), 
\eeq
so that it
is a critical point of the functional $\tilde J_n$ defined in 
$\H_{\eps_{n},x_{n}}$ by
\bdm\label{defJn}
\tilde J_{n}(u)=
\frac{1}{2}\|\nabla u\|^{2}_{2}+\frac{1}2\intr V(x_{n}+\eps_{n}x) u^{2}dx
-\int_{\R^{N}}G(x_{n}+\eps_{n}x,u(x))dx.
\edm
Notice that, by a simple change of scale and from \eqref{stimato},
it is possible to verify that
\begin{equation}\label{change}
\tilde J_{n}(\phi_{n})=\eps_{n}^{-N}J_{\eps_{n}}
(u_{\eps_{n}}), \qquad \|\phi_{n}\|^{2}_{H^{1}}\leq L.
\end{equation}
As in Lemmas \ref{ceramicomp} \ref{tresei} and from
elliptic regularity estimates, it results
that $\phi_n$  converges $C^2$ on compact sets to a function 
$\phi\in H^{1}$, which, by \eqref{contro} must be nontrivial.
Then, $\phi$ is a solution of
\bdm
-\Delta\phi+V(x_{0})\phi=\ov{ f}(x_{0},\phi(x)),
\edm
as $\chi(x_{0})=0$. 
This is the Euler equation of the functional
\bdm
\ov{I}_{x_0}(u)=\frac{1}{2}\|\nabla u\|_2^{2}+\frac{V(x_{0})}2
 \| u\|_{2}^{2}
-\int_{\R^{N}} \ov{F}(x_{0},u(x))dx.
\edm
On the other hand, conditions on $G$ allow us
to follow the same arguments of Lemma 2.2 in \cite{defe}
to deduce that 
\beq\label{claim}
\liminf_{n\to \infty}\tilde J_{n}(\phi_{n}) \geq \ov{I}_{x_0}(\phi).
\eeq
Indeed, consider the function
$$
h_{n}= \frac12\left[|\nabla \phi_{n}|^{2}+V(x_{n}+\eps_{n}x)
\phi_{n}^{2}\right] - G(x_{n}+\eps_{n}x,\phi_{n}(x)).
$$
Choosing $R>0$ sufficiently large, from the $C^{1}$ convergence of
$\phi_{n}$ over compacts, and since $\phi$  belongs
to $\hsob$ we have, for every $\delta>0$ fixed,
$$
\lim_{n\to\infty}\int_{B(0,R)}h_{n}\geq \ov{I}_{x_{0}}(\phi)-\delta.
$$
Moreover, taking $\eta_{R}$  a smooth cut-off function such that 
$\eta_{R}=0$ on $B(0,R-1)$ and $\eta_{R}= 1$ on $\RN\setminus B(0,R)$, 
and using as test function in \eqref{succ} $w=\eta_{R}\phi_{n}$, 
it is possible to obtain
$$
\liminf_{n\to\infty}\int_{\RN\setminus B(0,R)}h_{n}\geq -\delta,
$$
yielding \eqref{claim}.
Since $\phi$ is a critical point
of $\ov{I}_{x_0}$ and 
as the nonlinearity $\ov{f}(x_{0},t)/t$ 
is nondecreasing with respect to $t$,  
 we have 
\begin{equation}\label{maxiJ}
\ov{I}_{x_0}(\phi)=\max_{t\geq0}\ov{I}_{x_0}(t\phi).
\end{equation}
Moreover, it holds $F(x,t)\geq \overline{F}(x,t)$, for every 
$x\in \R^{N}$ and for every $t\in \R$,
so that,  Proposition 3.11 in \cite{ra}
together with \eqref{maxiJ}, implies that
\begin{equation}\label{diseq}
\ov{I}_{x_0}(\phi)=\max_{t\geq 0}\ov{I}_{x_{0}}(t\phi)
\geq \inf_{u\in H^{1}}\sup_{t\geq0} \ov{I}_{x_{0}}(tu)
\geq \inf_{u\in H^{1}}\sup_{t\geq0} I_{x_{0}}
(tu)=\Sigma(x_{0}).
\end{equation}
This inequality leads to an immediate contradiction 
in the case in which $x_{0}\not\in \Omega$, as it
would result $\Sigma(x_{0})=+\infty$ in this situation.
Otherwise we have that $x_{0}\in \Omega\cap B(z,r)$ and, 
assuming the first condition in \eqref{minivs} 
(the other case can be handled analogously) 
from Lemma \ref{decre} 
it follows that $F(x_{0},t)<F(z,t)$, moreover,
as $V(x_{0})\geq V(z)$ we have that $I_{x_{0}}(v)\geq I_{z}(v)$
for every $v\in \H$, which yields $\Sigma(x_0)>\Sigma(z)$, 
(for $\Sigma(z)$  defined  in \eqref{defsigmagen}). 
This, \eqref{leveleps}, \eqref{change}, 
\eqref{claim} and \eqref{diseq} yield
\begin{equation}
\label{cruc-step}
\Sigma(z)<\ov{I}_{x_0}(\phi)\leq  \liminf_{n\to\infty}\tilde J_{n}(\phi_{n})
\leq \liminf_{n\to +\infty}\eps_{n}^{-N}J_{\eps_{n}}(u_{\eps_{n}})
\leq \Sigma(z),
\end{equation}
which is a contradiction, proving \eqref{limeps}.

We are now ready to conclude the proof of the result. Let us fix $\delta>0$;
from \eqref{limeps} it follows that there exists $\eps_\delta>0$ such that
$0\leq u_{\eps}(x)<\delta$ 
for any $x\in\pa B(z,r)$ and $\eps\in(0,\eps_\delta)$.
It follows that $(u_{\eps}-\delta)^{+}= 0$  on
$\pa B(z,r)$ and hence we can choose
$$
\phi_\eps=(u_{\eps}-\delta)^{+}\chi_{\{|x-z|>r\}}\in H^1,
$$
as test functions in \rife{pbareps}. By multiplying and integrating
over $\R^{N}$, we obtain
$$
\int_{\R^{N}\diff B(z,r)}\left(\eps^{2}|\grad
(u_{\eps}-\delta)^{+}|^{2}+V(x)u_{\eps}(u_{\eps}-\delta)^{+}
-g(x,u_{\eps}(x))(u_{\eps}-\delta)^{+}\right)\\
=0.
$$
Having defined
\bdm
\Upsilon_\eps(x)=\begin{cases}
\dys V(x)-\frac{g(x,u_{\eps}(x))}{u_{\eps}}, & u_{\eps}(x)\neq 0
\\
0 & u_{\eps}(x)=0
\end{cases}
\edm
the preceding identity turns into
$$
 \int_{\R^{N}\diff B(z,r)}\big(\eps^{2}|\grad (u_{\eps}-\delta)^{+}|^{2}
+\Upsilon_\eps(x)|(u_{\eps}-\delta)^{+}|^{2}
+\Upsilon_\eps(x)\delta(u_{\eps}-\delta)^{+}\big) =0.
$$
By the definition of $g(x,t)$, it is easy to show that 
$\Upsilon_\eps(x)\geq 3\mu/4$
for all $x$ with $u_\eps(x)>0$, which implies that $(u_{\eps}(x)-\delta)^{+}=0$
 for every $x\not\in B(z,r)$ and every
$0<\eps<\eps_\delta$, namely the assertion.
\edim
\bos
The argument used in Proposition 
\ref{maximo}  has actually a stronger 
consequence, it implies that for every $\delta>0$
there exists $\eps_\delta>0$ such that
$$
\sup_{0<\eps<\eps_\delta}\,\sup_{x\in\R^{N}\setminus B(z,r)\cap \Omega}
u_{\eps}(x)<\delta.
$$
Indeed, using \eqref{limeps} and 
following the argument at the beginning of the proof 
we can prove
$$
\lim_{\eps\to 0}\sup_{\partial \Omega\cap B(z,r)}u_{\eps}=0
$$
arguing again   by contradiction and assuming that there exists 
a sequence $\{\eps_{n}\}$ converging to $0$ and a sequence 
$\{x_{n}\}\in \pa \Omega\cap B(z,r)$ such that, for some positive constant $\beta$,
\[
u_{\eps_n}(x_{n}) \geq \beta\qquad\text{for all $n\geq 1$.}
\]
As before, $x_{n}\to x_{0}\in \pa\Omega\cap B(z,r)$, and
$\phi_{n}(x)=u_{\eps_{n}}(x_{n}+\eps_{n}x)$ has a $C^{2}$ limit
$\phi$ satisfying
$$
-\Delta \phi+V(x_{0})\phi=g(x_{0},\phi).
$$
Since $x_{0}\in \pa\Omega$, $s(x_{0})V(x_{0})=1$, then,
applying Pohozaev identity to the solution
$\phi$, and recalling \eqref{fbar} we obtain
\begin{align*}
\frac{2}{2^{*}}\|\nabla \phi\|_{2}^{2}&=\|\phi\|_{2}^{2}V(x_{0})
(\chi(x_{0})-1)-\frac{\chi(x_{0})}{s^{2}(x_{0})}\intr \ln(1+s(x_{0})\phi^{2})+2(1-\chi(x_{0}))\intr \ov{F}(x_{0},\phi)
\\
&\leq -\frac{\chi(x_{0})}{s^{2}(x_{0})}\intr \ln(1+s(x_{0})\phi^{2})\leq 0
\end{align*}
giving the desired contradiction. 
\eos
{\bf Proof of Theorem \ref{mainvs}.}\\
By virtue of Proposition \ref{maximo}, taking into account the
definition of $G$,  $u_{\eps}$ turns out to be a
solution of \eqref{Eq:eps} for $\eps$ sufficiently
small. From elliptic regularity theory it follows that
$u_\eps$ is a positive $C^2$ function.  Let $\xi_\eps \in B(z,r)$ a local
maximum point of the function $u_{\eps}(x)$, then
$$
0\leq -\Delta u_{\eps}(\xi_\eps)=-V(\xi)u_{\eps}
(\xi_\eps)+f(\xi_{\eps},u_{\eps}(\xi_\epsilon))
\leq -V(\xi)u_{\eps}(\xi_{\eps})+u_{\eps}^{3}(\xi_{\eps})
$$
which implies that  there exists a positive constant $\sigma$,
independent on $\eps$, such that
\beq\label{stimainf}
u_{\eps}(\xi_\eps)\geq\sigma.
\eeq
Let us first prove conclusion {\it (ii)} of Theorem \ref{mainvs}
arguing by contradiction.
More precisely, consider $\eps_n\to 0$ and $x_n\in B(z,r)$
 a local maximum point of $u_{\eps_n}$.
Let $x_n\to x^*\in\overline{B}(z,r)$ and  consider the
sequence $\phi_n(x)=u_{\eps_n}(x_n+\eps_n x)$,
and its limit $\phi$, critical point of the autonomous 
functional $I_{x^*}$. Thanks to \eqref{stimainf},  
$\phi\neq 0$, implying that  $x^{*}\in \Omega$
Moreover, assuming the first alternative in
\eqref{minivs} (the other situation being similar)
 $s(x^*)>s_0$ and  $V(x^*)\geq V_0$, and
from Lemma  \ref{decre} we obtain that 
$F(x^{*},v)<F(z,v)$ so that
$I_{x^{*}}(v)> I_{z}(v)$ for
every $v\in H^{1}$, yielding the  inequality 
$$
I_{x^{*}}(\phi)\geq \inf_{u\in H^{1}}\max_{t>0}I_{x^{*}}(tu)
>\inf_{u\in H^{1}}\max_{t>0}I_{z}(tu)=\Sigma(z)
$$
which contradicts  \eqref{leveleps}, proving {\it(ii)}.
In order to prove conclusion {\it (i)} of Theorem \ref{mainvs}, 
assume by contradiction that there exist
a sequence $\{\eps_{n}\}$ converging to zero and  two local
maxima $x^{1}_{n}$, $x^{2}_{n}$ $\in \overline{B}(z,r)$, 
which both satisfy \eqref{stimainf}. We consider the sequence
$
\phi_{n}(x)=u_{\eps_{n}}(x^{1}_{n}+\eps_{n}x)
$
which is a critical point of the functional
$$
I^1_n(v)=\frac12\|\nabla v\|^2+\frac{1}{2}\intr V(x^{1}_{n}+\eps_{n}x)
v^2dx-\intr F(x_{n}^1+\eps_n x, v(x))dx
$$
with critical level 
\beq\label{uffa}
I^1_n(\phi_{n})=\eps_{n}^{-N}I_{\eps_{n}}(u_{\eps_{n}}).
\eeq
Arguing as before, we show that  $\phi_{n}$ converges in
the $C^{2}$ sense over compacts to a solution 
$\phi$ of \eqref{limity}
with $y=x^{1}$ and from conclusion {\it (ii)}
 $s(x^{1})=s_0$ and $V(x_{1})=V(0)$.
From \eqref{stimainf} we get that $\phi\neq 0$ and from
\cite{seta,sezo}
 we deduce that $\phi$ is a nonnegative, radially
symmetric function. Then, arguing as in the cubic case i.e.
$f(x,t)=t^{3}$,  $\phi$ has a local non-degenerate
maximum point, which, up to translations,  is located
in the origin. This facts and the $C^{2}$
convergence of $\phi_{n}$ imply that 
$x_n=(x^{1}_{n}-x^{2}_{n})/\eps_{n}\to
\infty$.  Then we can argue as in the proof of 
\eqref{claim} to get a contradiction.
Indeed, we consider the function
$$
h_n=\frac12|\nabla \phi_n|^2+\frac{1}2V(x^1_{n}+\eps_n x)
\phi_n^2-F(x^1_{n}+\eps_n x,\phi_n(x)),
$$
and observe that, thanks to the $C^2 $ convergence over compacts 
of $\phi_n$, for every $\delta$ we can choose $R>0$
\beq\label{dis1}
\lim_{n\to\infty}\int_{B(0,R)}h_n(x)dx  \geq I_{x^1}(\phi)-\delta.
\eeq
Moreover, as $x_n\to \infty$ we can fix $n_{0}$
sufficiently large such that $B(0,R)\cap B(x_n,R)=\emptyset$. 
On the other hand, the change of variable $y=x-x_n$ leads to
$$
\lim_{n\to\infty}\int_{B(x_n,R)}\!\!\!\! h_n (x)dx =\!\! \frac12\lim_{n\to\infty}
\int_{B(0,R)}\!\!\! |\nabla  \psi_n (y)|^2+V(x^2_{n}+\eps_n x) \psi_n^2(y)
 -2F(x^2_n+\eps_n y,{\psi}_n(y))dy
$$
where we put $\psi_{n}(y)=\phi_{n}(y+x_{n})$.
Reasoning as in \eqref{dis1} and taking into account that
$s(x^{1})=s(x^{2})=s_0$, we get
\beq
\label{dis2}
\lim_{n\to\infty}\int_{B(x_n,R)}h_n \geq
I_{x^{2}}(\psi)-\delta=I_{x^{1}}(\phi)-\delta.
\eeq
Then, arguing as in the proof of \eqref{claim} we get
$$
\liminf_{n\to\infty}\intr I^{1}_{n}(\phi_{n})\geq 2 \Sigma(x^{1})=2\Sigma(z),
$$
which is in contradiction with \eqref{uffa} and \eqref{leveleps}.
\vskip3pt
\noindent
In order to prove the exponential decay, notice that, by Proposition \ref{maximo},
$u_{\eps}$  decays  to zero at infinity, uniformly with respect to $\eps$. 
Hence we find $\rho>0$, $\Theta\in(0,\sqrt{\mu})$ and $\eps_0>0$ 
such that $u_\eps^2 \leq \mu-\Theta^2$,
for all $|x-x_\eps|>\eps\rho$ and $0<\eps<\eps_0$. Let us set
$$
\xi_\rho(x)=M_\rho e^{-\Theta(|x-x_\eps|/\eps-\rho)},\qquad
M_\rho=\sup_{(0,\eps_0)}\max_{|x|=\rho} (u_\eps),
$$
and introduce the set $\dys {\mathcal A}=\bigcup_{R>\rho}D_R,$
where, for any $R>\rho$,
$$
\quad D_R=\big\{\rho<|x|<R:\,\,
u_\eps(x)>\xi_\rho(x)\,\,\,\,
\text{for some $\eps\in(0,\eps_0)$}\big\}.
$$
Assume by contradiction that ${\mathcal A}\not=\emptyset$.
Then there exist $R_*>\rho$ and $\eps_*\in(0,\eps_0)$ with
$$
\eps^2\Delta(\xi_\rho-u_{\eps_*}) 
\leq\left[\Theta^2-\frac{2\eps\Theta}{|x-x_\eps|}\right]\xi_\rho
-u_{\eps_{*}}(\mu-u^{2}_{\eps_{*}})\leq 
\Theta^2 (\xi_{\rho}-u_{\eps_*} )<0,
$$
for $x$ in $D_R$ and for all $R\geq R_*$.
Hence, by the maximum principle, we get
$$
\xi_\rho-u_{\eps_*} \geq
\min\Big\{\min_{|x|=\rho}(\xi_\rho-u_{\eps_*} ),
\min_{|x|=R}(\xi_\rho-u_{\eps_*} )\Big\},
$$
in $D_R$ for all $R\geq R_*$. Letting $R\to\infty$ and
recalling the definition of $\xi_\rho$ yields
$$
\xi_\rho-u_{\eps_*}\geq \min
\Big\{\min_{|x|=\rho}(\xi_\rho-u_{\eps_*}),0\Big\}\geq 0,
\quad\text{in $\bigcup_{R\geq R_*}D_R$}.
$$
In turn, $u_{\eps_*}(x)\leq \xi_\rho(x)$ for
all $x$ in $\dys \bigcup_{R\geq R_*}D_R$, which yields a contradiction.
Whence ${\mathcal A}=\emptyset$, and the desired exponential decay follows.
\edim
{\bf Proof of Theorem \ref{mainsigma}.}\\
Theorem \ref{mainsigma} can be proved as Theorem \ref{mainvs}; 
indeed one can perform the penalization procedure around the minimum
point of $\Sigma$ introduced in \eqref{minisigma} and make the 
analogous calculation up to \eqref{cruc-step}, as it is readily seen that
hypothesis \eqref{minivs} is used to obtain condition \eqref{minisigma}.
The rest of the proof can be handled in the same way as in the proof of
Theorem \ref{mainvs}.
\edim
\section{Proof of Theorem \ref{mainnece}}
As a first step  to prove Theorem \ref{mainnece}, let us show the following Lemma.
\ble\label{pucci}
Let us suppose that $V,\,s\in C^{1}(\R^{N})$ satisfies \eqref{growths}.
If $z\in {\mathcal E}$, then
$$
\nabla \left(V(z)-\frac1{s(z)}\right)\|Q_{z}\|_{2}^{2}
+\intr \nabla\left(\frac1{s^{2}(z)}\ln(1+s(z)Q^{2}_{z}(x))\right)dx=0
$$
where $Q_{z}$ is the least energy solution of the autonomous Problem
\eqref{limity} for $y=z$.
\ele
\bdim
We will closely follow the argument in \cite{mopesq} (see also \cite{sesq}).
Let $z\in {\mathcal E}$, a sequence $\{\eps_{n}\}$ converging to zero 
and $u_{\eps_{n}}$ a solution of Problem \eqref{Eq:eps}, for
$\eps=\eps_{n}$, as in Definition \ref{defe}.
Let us define $\vfi_{n}(x)=u_{\eps_{n}}(z+\eps_{n}x)$ 
and apply the Pucci--Serrin identity \cite[Proposition 1]{puse} 
with the lagrangian function
${\mathcal L}:\R^{N}\times \R\times \R^{N}\to \R$ defined by
$$
{\mathcal L}(x,t,\xi)=\frac12|\xi|^{2}+
V(z+\eps_{n}x)\frac{t^{2}}{2}-F(z+\eps_{n}x,t),
$$
obtaining
\begin{align*}
\sum\limits^N_{i,\ell=1} \intr\partial_{x_{i}} {\boldsymbol h}^\ell  
\partial_{x_{i}} \vfi_n \, \partial_{x_{\ell}} \vfi_n  dx&
=\intr {\rm div} {\boldsymbol h}\, {\mathcal L}(x,\vfi_n,\nabla \vfi_n)
+
\frac{\eps_n}2\intr {\boldsymbol h} \cdot  
\nabla_{x} V(z+\eps_n x)\varphi_n^2 
\\
&+\eps_{n}\intr {\boldsymbol h}\cdot \nabla_{x} F(z+\eps_{n}x,\vfi_{n}(x))
\end{align*}
for all ${\boldsymbol h} \in C^1_{\rm c}\left(\R^{N}, \R^{N} \right)$.
Let us choose, for any $\lambda>0$,
$$
{\boldsymbol h}_j:\R^N\to\R^N,\qquad
{\boldsymbol h}_j^\ell(x)=
\begin{cases}
\Upsilon(\lambda x) & \text{if $\ell=j$}, \\
0 & \text{if $\ell\neq j$},
\end{cases}\qquad \ell=1,\dots, N
$$
with $\Upsilon\in C^1_{{\rm c}} (\R^{N})$,
$\Upsilon(x)=1$ if $|x|\leq 1$ and $\Upsilon(x)=0$
if $|x|\geq 2$. Then, for $j=1,\dots,N$,
\begin{align*}
& \sum\limits^N_{i=1}
\intr \lambda \partial_{x_{i}} \Upsilon(\lambda x) 
\partial_{x_{i}} \vfi_n \, \partial_{x_{j}} \vfi_n
=\intr \lambda \partial_{x_{j}} \Upsilon(\lambda x) {\mathcal L}(x,\vfi_n,
\nabla \vfi_n) \\
&+\frac{\eps_n}2\intr  
\Upsilon(\lambda x)[\partial_{x_{j}} V(z+\eps_n x)
\varphi_{n}^{2}-2\partial_{x_{j}}F(z+\eps_{n}x,\varphi_{n}(x))].
\end{align*}
By the arbitrariness of $\lambda>0$, letting
$\lambda \to 0$ and keeping $j$ fixed, we obtain
\[
\intr [\partial_{x_{j}} V(z+\eps_n x)\varphi_n^2-2
\partial_{x_{j}} F(z+\eps_n x,\varphi_n(x))]dx=0
\qquad j=1,\dots,N.
\]
By assumption \eqref{growths}, there exists a positive constant 
$\beta_1$ such that, for all $x\in\R^N$ and $j\geq 1$, we get
$|\nabla V(z+\eps_n x)|\leq \beta_1 e^{\gamma\eps_n |x|}$
and $|\nabla s(z+\eps_n x)|\leq \beta_1 e^{\gamma\eps_n |x|}$,
so that, invoking the uniform exponential decay of $\varphi_n$,
letting $n\to \infty$ in the above identity, and
recalling that $\varphi\to Q_{z}$, the least
energy solution of \rife{limity} 
for $y=z$, we find
\begin{equation}\label{part-fin}
\intr [\partial_{x_{j}} V(x){\big |}_{z}Q_z^2-2
\partial_{x_{j}} F(x,Q_z(x)){\big |}_{z}]dx=0
\qquad j=1,\dots,N,
\end{equation}
giving the conclusion.
\edim
The following result will be useful in studying the function $\Sigma$.
\begin{lemma} \label{theta}
Assume \eqref{s0} and that $V,\,s \in C^{1}(\R^{N})$. Then 
The function ${\mathcal G}: H^{1}\times \R^{N}\to \R$
defined by
$$
\mathcal{G}(u,y):= \|\nabla u\|_{2}^{2}+ \left(V(y)-\frac1{s(y)}\right) 
\|u\|^{2}_{2}
$$  
is continuous in $y$, for any $u \in H^1(\RN)$.
Moreover, if $\; \mathcal{G}(u,y) < 0$, then there exists  a unique 
$\theta(u,y) > 0$ such that $\theta(u,y)  u \in \mathcal{N}_y$, where 
$\mathcal{N}_y$ is defined in \eqref{ny}. Finally,
the map $\theta$ is continuous on $H^{1}(\R^{N})\times \R^{N}$
\end{lemma}
\bdim
Given $u$ such that  $\; \mathcal{G}(u,y) < 0$ , we define the function 
$g: [0, \infty)\times H^{1}(\R^{N})\times \R \to \R$ as follows 
$$
g(\tau,u,y):=\begin{cases}
\dys \|\nabla u\|_{2}^{2}+ V(y) \|u\|_{2}^{2}\; & \text{ for } \tau=0,
\medskip\\
\dys \frac{1}{\tau^{2}}\langle ( I_y)^\prime(\tau u), 
\tau u \rangle  &\text { for } \tau > 0.
\end{cases} 
$$
Condition \eqref{s0} and the regularity properties of the functions
$V$ and $s$ imply that $g$ is continuous and
Lebesgue Dominate Convergence Theorem yields
$$
\lim_{\tau \to +\infty} g(\tau)= \|\nabla u\|_{2}^{2}+ 
 \left(V(y)-\frac1{s(y)}\right)  \|u\|_{2}^{2}< 0.\;
$$
Since $g$ is a continuous function, there exists  
$\theta(u,y) > 0$ such that $g( \theta(u,y)) =0$, that is
\[
\langle (I_y)^\prime(\theta(u,y) u), \theta(u,y) u \rangle =0,
\]
i.e. $\theta(u,y) u \in \mathcal{N}_y$. The uniqueness of 
$\theta(u,y)$ follows from the fact that 
$f(y,u)=\partial_{u}F(y,u)$ satisfies
$f(u,y)/u$ is nondecreasing with respect to $u$.
\\
The continuity of the $\theta$ can be deduced from the
Implicit function Theorem. 
\edim
In order to prove Theorem \ref{mainnece} let us first show the 
regularity properties of the function $\Sigma$.
\bpr\label{sigmadiff}
Assume \eqref{s0} and that $V,\,s\in C^{1}(\R^{N})$. Then,
the function $\Sigma$ is of class $C^{1}(\Omega)$ and its gradient is given by
\begin{equation}\label{gradsigma}
\nabla \Sigma(y)=\nabla_y \left(V(y)-\frac{1}{s(y)} \right)
\| Q_y\|_{2}^{2}
+\intr \nabla_y \left( \frac{1}{s(y)^{2}}\ln(1+s(y)Q_y^2(x))\right)\;dx \,,
\end{equation}
where $Q_{y}$ is the least energy solution of \eqref{limity}.
\epr
\bdim
In order to compute the directional derivative of the function 
$\Sigma$ with respect to a unitary vector $\eta \in \RN$,
let $\rho =y+\tau\eta$, and as $\tau\to 0$ then $\rho\to y$.

Since $\mathcal{G}(u,y)$ is continuous in $y$ and 
$\mathcal{G}(Q_y,y) < 0$, by Lemma \ref{theta}, there exists a 
$\delta > 0$ such that, for $|\tau|=\| \rho - y\| < \delta$, 
then $\mathcal{G}(Q_{y},\rho) < 0$. By Lemma \ref{theta}, there 
exists $\theta (Q_{y},\rho)=\theta (y,\rho)>0$ such that 
$\theta (y,\rho) Q_y \in \mathcal{N}_{\rho}$.
Using the Mean Value Theorem and the definition \eqref{defsigmagen}, 
we have
\begin{equation}\label{differC1}
\Sigma(\rho)-\Sigma(y) \leq 
I_\rho(\theta (y,\rho) Q_y) - I_y(\theta(y,y)Q_y)
= \tau\eta\cdot\nabla_\xi I_\xi (\theta(y,\xi)Q_y) 
|_{\xi \in [y,y+ \tau\eta]}
\end{equation}
Computing $\nabla_\xi I_\xi (\theta(y,\xi)Q_y)$, we obtain
\begin{eqnarray*}
\nabla_\xi I_\xi (\theta(y,\xi)Q_y)&=&\frac{\theta(y,\xi)^2 }2 
\nabla_\xi\left[ V(\xi) - \frac{1}{s(\xi)}
\right]
\| Q_y\|_{2}^2
\\
&&+  \frac12 \nabla_\xi \left(\frac{1}{s^2(\xi)}\right) 
\intr \ln(1+ s(\xi) \theta^2(y,\xi) Q_y^2)dx\\
&&+ \frac12 \frac{1}{s^2(\xi)} \intr  \nabla_\xi 
\ln(1+ s(\xi) \theta^2(y,\xi) Q_y^2)dx\\
&&+ \nabla_\xi \theta(y,\xi)
\theta(y,\xi)\left[ \|\nabla Q_y\|_{2}^{2}+\Big(V(\xi)-\frac1{s(\xi)}\Big)
\|Q_y\|_{2}^{2}\right] \;. 
\end{eqnarray*}
Since $\theta(y,\xi)Q_{y}\in {\mathcal N}_{\xi}$ it results
$\nabla_\xi \theta(y,\xi)\langle 
(I_\xi)^\prime (\theta(y,\xi)Q_y), \theta(y,\xi) Q_y \rangle/
\theta(y,\xi)=0$ and, on the other hand
\[
\langle 
(I_\xi)^\prime (\theta(y,\xi)Q_y), \theta(y,\xi) Q_y \rangle
=
\|\theta\nabla Q_{y}\|_{2}^{2}+\|\theta V(\xi) Q_{y}\|_{2}^{2}-\intr
\frac{\theta^{4}Q_{y}^{4}}{1+\theta^{2}s(\xi)Q_{y}^{2}}
\]
and substituting above, we get 
\begin{eqnarray*}
\nabla_\xi I_\xi (\theta(y,\xi)Q_y)&\!\!=& \frac{\theta(y,\xi)^2}2  
\nabla_\xi\left(V(\xi)-\frac{1}{s(\xi)}\right) 
\| Q_y\|_{2}^2\\
&\quad +& \frac12 \nabla_\xi \left(\frac{1}{s^2(\xi)}\right) 
\intr \ln(1+s(\xi) \theta^2(y,\xi) Q_y^2)\\
 &\quad+&\frac1{2s^{2}(\xi)}\intr \nabla_{\xi}\ln(1+s(\xi)\theta^{2}(y,\xi)
Q_{y}^{2})
\\
&\quad+&\frac{\nabla_\xi \theta(y,\xi)}{\theta(y,\xi)} 
\intr\frac{\theta^{4}(y,\xi)Q_{y}^{4}}{1+s(\xi)\theta^{2}(y,\xi)Q_{y}^{2}}-\frac{\theta^{2}(y,\xi)Q_{y}^{2}}{s(\xi)}
 \;. 
\end{eqnarray*}
Since
\begin{eqnarray*}
\frac1{2s^{2}(\xi)}\intr \nabla_{\xi}\ln(1+s(\xi)\theta^{2}(y,\xi)
Q_{y}^{2}) &\!\!=&
\frac{\nabla s(\xi)}{2s^{2}(\xi)}\intr \frac{\theta^{2}(y,\xi)Q_{y}^{2}}{1+s(\xi)\theta^{2}(y,\xi)Q_{y}^{2}}\\
&\quad +&\frac{\nabla_\xi \theta(y,\xi)}{\theta(y,\xi)}
\intr \frac{\theta^{2}(y,\xi)Q_{y}^{2}}{s(\xi)(1+s(\xi)\theta^{2}(y,\xi)Q_{y}^{2})}
\end{eqnarray*}
substituting above, we get 
\begin{eqnarray*}
\nabla_\xi I_\xi (\theta(y,\xi)Q_y)&=& \frac{\theta(y,\xi)^2}2  
\nabla_\xi\left(V(\xi)-\frac{1}{s(\xi)}\right) 
\| Q_y\|_{2}^2\\
&& + \frac12 \nabla_\xi \left(\frac{1}{s^2(\xi)}\right) 
\intr \ln(1+s(\xi) \theta(y,\xi)^2 Q_y^2)\\
&&+ \frac12 \frac{\nabla_\xi s(\xi)}{s^2(\xi)} \intr 
\frac{\theta(y,\xi)^2 Q_y^2 }{1 + s(\xi) \theta(y,\xi)^2 Q_y^2} \;. 
\end{eqnarray*}
Using (\ref{differC1})  we obtain
\[
\limsup_{\tau \to 0^+} \frac{\Sigma(y+\tau\eta)-\Sigma(y)}{\tau} 
\leq\limsup_{\tau \to 0^+}\;\eta \cdot \nabla_\xi I_\xi 
(\theta(y,\xi)Q_y) |_{\xi \in [y,y+\tau\eta]} \;.
\]
Taking into consideration that $\tau \to 0^+$ implies that 
$\xi\to y$, applying Lemma \ref{theta}, and observing that
$\theta(y,y)Q_y=Q_y\in\mathcal{N}_y$, 
we obtain
\begin{eqnarray*}
\limsup_{\tau \to 0^+} \frac{\Sigma(y+\tau\eta)-\Sigma(y)}{\tau}
	& \leq & \eta\cdot\nabla_\xi I_\xi (\theta(y,\xi)Q_y) |_{\xi=y}
	=  \frac12\eta\cdot\nabla_y \left(V(y)-\frac{1}{s(y)} 
\right) \|Q_y\|_{2}^2 \\
& &  +\frac12\intr \eta\cdot\nabla_y \left( 
\frac{1}{s(y)^{2}}\ln(1+s(y)Q_y^2(x))\right)\;dx\;.
\end{eqnarray*}

On the other hand, using the Mean Value Theorem and the definition
\eqref{defsigmagen}, we have
\begin{equation}\label{differC}
\Sigma(\rho)-\Sigma(y) \geq I_\rho(\theta(\rho,\rho)Q_\rho) 
- I_y(\theta(y,\rho)Q_\rho)
= \tau \eta \cdot \nabla_\xi I_\xi (\theta(\xi,\rho)Q_\rho) 
|_{\xi\in[y,y+\tau\eta]}\;.
\end{equation}
Performing a similar argument and observing that $\rho=y+\tau\eta\to y$, 
as $\tau \to 0^+$, yields
\begin{eqnarray*} 
\liminf_{\tau \to 0^+} \frac{\Sigma(y+\tau\eta)-\Sigma(y)}{\tau}
& \geq & \eta\cdot\nabla_\xi I_\xi(\theta(y,\xi)Q_y) |_{\xi=y}
=   +\frac12\eta\cdot\nabla_y 	\left(V(y)-\frac{1}{s(y)} \right) 
\|Q_y\|_{2}^2 \\
& & +\frac12\intr \eta\cdot\nabla_y \left( 
\frac{1}{s(y)^{2}}\ln(1+s(y)Q_y^2(x))\right)\;dx\;.
\end{eqnarray*}
The two inequalities give 
\[
\left(\frac{\partial \Sigma}{\partial \eta}\right)^+\!\!\!\!\!\!(y) 
=   +\frac12\eta\cdot\nabla_y 	\left(V(y)-\frac{1}{s(y)} \right) 
\|Q_y\|_{2}^2 
 +\frac12\intr \eta\cdot\nabla_y \left( 
\frac{1}{s(y)^{2}}\ln(1+s(y)Q_y^2(x))\right)\;dx\;.
\]
Analogously, we can prove that the same indentity holds for
\((\partial \Sigma/\partial \eta)^-(y)\), giving the directional
derivative of \(\Sigma\) along a vector $\eta$, showing 
\eqref{gradsigma}. The continuity of the gradient
follows from the regularity properties of $V,\,s$, \eqref{growths}
and the exponential decay of $Q_{y}$.
\edim

{\bf Proof of Theorem \ref{mainnece}.}\\
The regularity property of $\Sigma$ are proved in Proposition \ref{sigmadiff}.
Let now $z\in {\mathcal E}$, then Lemma \ref{pucci} and Proposition 
\ref{sigmadiff}  imply that $z$ is a critical point of $\Sigma$.
In order to show that $\nabla V$ and $\nabla s$ are linearly dependent,
let us compute the following partial derivative
\begin{eqnarray*}
\intr \partial_{j}\left(\frac1{s^{2}(z)}
\ln(1+s(z)Q_{z}^{2}(x))\right)dx &=&
\partial_{j}\left(\frac1{s^{2}(z)}\right)\intr 
\ln(1+s(z)Q_{z}^{2}(x))dx
\\
&&+\frac{\partial_{j}s(z)}{s^{2}(z)}
\intr\frac{Q^{2}_{z}(x)}{1+s(z)Q^{2}_{z}(x)}dx
\\
&&=\frac2{s(z)}\partial_{j}\left(\frac1{s(z)}\right)
\intr\ln(1+s(z)Q_{z}^{2}(x))dx
\\
&&-\partial_{j}\left(\frac1{s(z)}\right)
\intr\frac{Q^{2}_{z}(x)}{1+s(z)Q^{2}_{z}(x)}dx. 
\end{eqnarray*}
Therefore, every $z\in{\mathcal E}$ has to satisfy
\[
\partial_{j} V(z)\|Q_{z}\|_{2}^{2}=
- \frac{\partial_{j} s(z)}{s^{3}(z)}
\intr \left[s(z) Q_{z}^{2}(x)-2\ln(1+s(z)Q_{z}^{2}(x))
+\frac{s(z)Q^{2}_{z}(x)}{1+s(z)Q^{2}_{z}(x)}\right] dx,
\]
for all $j=1,\dots,N$,
showing the linear dependence of $\nabla V$ and $\nabla s$.
The proof is complete once one takes into account that the 
function $h(t)=t-2\ln(1+t)+t/(1+t)$ is positive for $t>0$.
\edim

\bos\label{fine}
We can precise Remark \ref{quasifine},
in the sense that if $z$ is not a critical point of $V$ and $s$,
there exists at least a $j_{0}$ such that 
$\partial_{j_{0}} V(z)\neq 0$ and $\partial_{j_{0}}s(z)\neq 0$, 
nevertheless $\partial_{j_{0}}V(z)$
and $\partial_{j_{0}}s(z)$ have to  satify \eqref{uffa}.
\eos



\begin{thebibliography}{99}

\bibitem{ambaci} A. Ambrosetti, M. Badiale, S. Cingolani,
{\it Semiclassical states of nonlinear Schr\"odinger equations}, 
Arch. Rational Mech. Anal. {\bf 140} (1997) 285-300.

\bibitem{amfema} A. Ambrosetti, V. Felli, A. Malchiodi,
{\it Ground states of nonlinear Schr\"odinger equations 
with potentials vanishing at infinity}, J. Eur. Math. Soc. 
(JEMS) {\bf 7} (2005) 117-144.

\bibitem{babefo} P. Bartolo, V. Benci, D. Fortunato, 
{\it Abstract critical point theorems and applications 
to some nonlinear problems with strong resonance at 
infinity}, Nonlinear Anal. TMA {\bf 7} (1983) 981-1012.

\bibitem{beje} A. Beyon, L. Jeanjean,
{\it Standing waves for nonlinear Schr\"odinger equations 
with a general nonlinearity}, Arch. Ration. Mech. Anal. 
{\bf 185} (2007) 185-200.

\bibitem{beje2} A. Beyon, L. Jeanjean,
{\it erratum to: "Standing waves for nonlinear Schr\"odinger equations 
with a general nonlinearity"}, Arch. Ration. Mech. Anal. 
{\bf 190} (2008) 549-551.

\bibitem{bova} D. Bonheure, J. Van Schaftingen,
{\it Bound state solutions for a class of nonlinear 
Schr\"odinger equations}, Rev. Mat. Iberoam. 
{\bf 24} (2008) 297-351.

\bibitem{cerami} G. Cerami, 
{\it Un criterio di esistenza per i punti critici su variet\`a
illimitate}, Rend. Accad. Sc. Lett. Inst. Lombardo 
{\bf 112} (1978) 332-336.

\bibitem{defe} M. Del Pino, P. Felmer,
{\it Local mountain passes for semi-linear elliptic problems
in unbounded domains}, Calc. Var. Partial Differential Equations 
{\bf 4} (1996) 121-137.

\bibitem{defe2} M. Del Pino, P. Felmer,
{\it Semi-classical states for nonlinear Schr\"odinger equations},
J. Funct. Anal. {\bf 149} (1997) 245-265.

\bibitem{flwe} A. Floer, A. Weinstein,
{\it Nonspreading wave packets for the cubic
Schr\"odinger equation with a bounded potential},
J. Funct. Anal. {\bf 69} (1986) 397-408.

\bibitem{gr}  M. Grossi, 
{\it  Some results on a class of nonlinear Schr\"odinger equations}, 
Math. Zeit. {\bf 235} (2000) 687-705.

\bibitem{jeta1} L. Jeanjean, K. Tanaka, 
{\it A remark on least energy solutions in $R^N$}, 
Proc. Amer. Math. Soc. {\bf 131} (2003) 2399-2408.

\bibitem{jeta2} L. Jeanjean, K. Tanaka, 
{\it Singularly perturbed elliptic problems with super linear 
or asymptotically linear nonlinearities}, Calc. Var. 
Partial Differential Equations {\bf 21} (2004) 287-318.

\bibitem{li} Y.Y. Li, 
{\it  On a singularly perturbed elliptic equation}, 
Adv. Differential Equations {\bf 2} (1997) 955-980.

\bibitem{plions} P.L. Lions,
{\it The concentration-compactness principle in the calculus of
variations. The locally compact case}, Ann. Inst. H. Poincar\'e, 
Anal. Non Lin\'eaire {\bf 1} (1984) 223-283.

\bibitem{mamopesat} L.A. Maia, E. Montefusco, B. Pellacci, 
{\it Weakly coupled nonlinear Schr\"odinger systems: the 
saturation effect}, Calc. Var. Partial Differential Equations 
{\bf 46} (2013) 325-351.

\bibitem{mopesq} E. Montefusco, B. Pellacci, M. Squassina,
{\it Semiclassical states for weakly coupled nonlinear 
Schr\"odinger systems}, 
J. Eur. Math. Soc. {\bf 10} (2008) 47-71.

\bibitem{oh1}  Y.G. Oh, 
{\it Existence of semiclassical bound states of nonlinear
Schr\"odinger equations with potentials of the class $(V)_{a}$}, 
Comm. Partial Differential Equations {\bf 13} (1988) 1499-1519.

\bibitem{oh2} Y.G.  Oh, 
{\it On positive multi-lump bound states of nonlinear 
Schr\"odinger equations under multiple well potential}, 
Comm. Math. Phys. {\bf 131} (1990) 223-253.

\bibitem{pi}  A. Pistoia, 
{\it  Multi-peak solutions for a class of nonlinear 
Schr\"odinger equations}, NoDEA Nonlinear Diff. Eq. 
Appl. {\bf 9} (2002) 69-91.

\bibitem{puse}  P. Pucci, J. Serrin, 
{\it  A general variational identity}, Indiana Univ. 
Math. J. {\bf 35} (1986) 681-703.

\bibitem{ra}  P.H. Rabinowitz, 
{\it  On a class of nonlinear Schr\"odinger equations}, 
Z. Angew. Math. Phys. {\bf 43} (1992) 270-291.

\bibitem{sesq}  S. Secchi, M. Squassina, 
{\it On the location of concentration points for singularly 
perturbed elliptic equations}, Adv. Differential Equations
{\bf 9} (2004) 221-239.

\bibitem{seta}  J. Serrin, M. Tang, 
{\it Uniqueness of ground states for quasilinear elliptic equations}, 
Indiana Univ. Math. J. {\bf 49} (2000) 897-923.

\bibitem{sezo}  J. Serrin, H. Zou, 
{\it Simmetry of ground states of quasilinear elliptic equations}, 
Arch. Ration. Mech. Anal. {\bf 148} (1999) 265-290.

\bibitem{stzo} C.A. Stuart, H. Zou,
{\it Applying the mountain pass theorem to an asymptotically 
linear elliptic equation on $\R^{N}$}, 
Comm. Partial Differential Equations {\bf 24} (1999) 1731-1758.

\bibitem{wangzeng} X. Wang, B. Zeng,
{\it On concentration of positive bound stated of
nonlinear Schr\"odinger equations with competing 
potential functions},
SIAM J. Math. Anal. {\bf 28} (1997) 633-655.
\bigskip

Liliane de Almeida Maia,\\
Departmento de Matem\'atica,\\
Uni\-ver\-si\-da\-de de Bras\'\i{}lia,\\ 
70910-900 Brasilia, Brazil.\\
E-mail address: {\tt lilimaia@unb.br}
\medskip

Eugenio Montefusco,\\
Dipartimento di Matematica, \\
{\it Sa\-pien\-za} Universit\`a  di Roma,\\ 
p.le Aldo Moro 5, 00185 Roma, Italy.\\
E-mail address: {\tt eugenio.montefusco@uniroma1.it}
\medskip

Benedetta Pellacci,\\ 
DiST Dipartimento di Scienze e Tecnologie,\\ 
Universit\`a degli Studi di Napoli {\it Parthenope},\\
Centro Direzionale Isola C4, 80143 Napoli, Italy.\\ 
E-mail address: {\tt pellacci@uniparthenope.it}

\end{thebibliography}
\end{document}